\documentclass{amsart}

\usepackage{amssymb}
\usepackage{amsmath}
\usepackage{amsthm}
\usepackage{amsfonts}

\usepackage{a4wide}

\newtheorem{proposition}{Proposition}
\newtheorem{theorem}{Theorem}
\newtheorem{corollary}{Corollary}

\theoremstyle{definition}
\newtheorem{definition}{Definition}
\newtheorem{example}{Example}
\newtheorem{experiment}{Experimental result}
\newtheorem{remark}{Remark}

\numberwithin{equation}{section}

\begin{document}


\title[Bicyclic commutator quotients]
{Bicyclic commutator quotients with \\ one non-elementary component}

\author{D. C. Mayer}
\address{Naglergasse 53\\8010 Graz\\Austria}
\email{algebraic.number.theory@algebra.at}
\urladdr{http://www.algebra.at}


\thanks{Research supported by the Austrian Science Fund (FWF): projects J0497-PHY and P26008-N25}

\subjclass[2010]{Primary 11R37, 11R32, 11R11, 11R20, 11R29, 11Y40;
20D15, 20E18, 20E22, 20F05, 20F12, 20F14}

\keywords{Hilbert \(3\)-class field tower, maximal unramified pro-\(3\) extension,
unramified cyclic cubic extensions, Galois action, imaginary quadratic fields,
bicyclic \(3\)-class group, punctured capitulation types, statistics;
pro-\(3\) groups, finite \(3\)-groups, generator rank, relation rank, Schur \(\sigma\)-groups,
low index normal subgroups, kernels of Artin transfers, abelian quotient invariants,
\(p\)-group generation algorithm, descendant trees, antitony principle}

\date{Sunday, August 22, 2021}


\begin{abstract}
For any number field \(K\) with non-elementary \(3\)-class group
\(\mathrm{Cl}_3(K)\simeq C_{3^e}\times C_3\), \(e\ge 2\),
the \textit{punctured capitulation type} \(\varkappa(K)\) of \(K\)
in its unramified cyclic cubic extensions \(L_i\), \(1\le i\le 4\),
is an \textit{orbit under the action of} \(S_3\times S_3\).
By means of Artin's reciprocity law,
the arithmetical invariant \(\varkappa(K)\) is translated
to the \textit{punctured transfer kernel type} \(\varkappa(G_2)\)
of the automorphism group \(G_2=\mathrm{Gal}(\mathrm{F}_3^2(K)/K)\)
of the second Hilbert \(3\)-class field of \(K\).
A classification of finite \(3\)-groups \(G\)
with low order and bicyclic commutator quotient
\(G/G^\prime\simeq C_{3^e}\times C_3\), \(2\le e\le 6\),
according to the algebraic invariant \(\varkappa(G)\),
admits conclusions concerning the length of the
Hilbert \(3\)-class field tower \(\mathrm{F}_3^\infty(K)\) of 
imaginary quadratic number fields \(K\).
\end{abstract}

\maketitle


\section{Introduction}
\label{s:Intro}

\noindent
An indispensable tool
for the investigation of the unramified Hilbert \(3\)-class field tower \(\mathrm{F}_3^\infty(K)\)
of an arbitrary number field \(K/\mathbb{Q}\)
with bicyclic \(3\)-class group \(\mathrm{Cl}_3(K)\simeq C_{3^e}\times C_3\), \(e\ge 2\),
is the \textit{punctured capitulation type}
\(\varkappa(K)=(\ker(\tau_i))_{1\le i\le 4}\)
of \(K\) in its four unramified cyclic cubic extensions \(L_1,\ldots,L_3;L_4\).
The puncture at \(L_4\) is motivated by the special role of \(L_4\)
as common subfield of all four unramified extensions
of degree nine of \(K\).
Here we denote by \(\tau_i:\,\mathrm{Cl}_3(K)\to\mathrm{Cl}_3(L_i)\),
\(\mathfrak{a}\mathcal{P}_K\mapsto(\mathfrak{a}\mathcal{O}_{L_i})\mathcal{P}_{L_i}\),
the \textit{transfers} (extension homomorphisms)
of \(3\)-classes of \(K\) into \(L_i\) for \(1\le i\le 4\).
We expand these ideas exemplarily for \(e=2\).
With minor modifications, however, they may be adopted for any \(e\ge 3\).

By means of Artin's reciprocity law
\cite{Ar1927,Ar1929},
we translate the arithmetical invariant \(\varkappa(K)\)
to the \textit{punctured transfer kernel type} \(\varkappa(G_2)\)
of the automorphism group \(G_2=\mathrm{Gal}(\mathrm{F}_3^2(K)/K)\)
of the second Hilbert \(3\)-class field of \(K\).
Based on the lattice of normal subgroups between
a pro-\(3\) group \(G\) with \(G/G^\prime\simeq C_9\times C_3\)
and its commutator subgroup \(G^\prime\) in \S\
\ref{ss:TopEndSbgrp},
we define the group theoretic \textit{Artin transfers}
\(T_i:\,G/G^\prime\to H_i/H_i^\prime\)
from \(G\) to maximal subgroups \(H_i\) with \((G:H_i)=3\),
corresponding to the arithmetical transfers \(\tau_i\), in \S\
\ref{ss:Transfers}.

We explain in \S\
\ref{ss:Orbits}
why only the \textit{orbit under the action of} \(S_3\times S_3\)
of the punctured transfer kernel type
\(\varkappa(G)=(\ker(T_i))_{1\le i\le 4}\)
is an \textit{invariant} of \(G\),
but not an individual orbit representative.

We conclude these preliminaries in \S\
\ref{ss:pTKT}
with an overview of all \(52\) combinatorially possible orbits,
emphasizing those which can actually be realized as \(\varkappa(G)\) by a \(3\)-group \(G\).

In \S\
\ref{s:NineByThree}
we devote our attention to a collection of algebraic invariants of the smallest
metabelian \(3\)-groups \(G\) with commutator quotient \(G/G^\prime\simeq C_9\times C_3\)
and coclass \(2\le\mathrm{cc}(G)\le 3\),
in order to establish the required inventory of
Galois groups \(\mathrm{Gal}(\mathrm{F}_3^2(K)/K)\)
for imaginary quadratic fields \(K\)
in \S\S\
\ref{s:ImagQuadratic}
and
\ref{s:HigherOrder},
supplemented by \S\
\ref{s:TwentysevenByThree}
and periodic Schur \(\sigma\)-groups in \S\S\
\ref{s:Metabelian}
and
\ref{s:PeriodicityAndLimits}.

In \S\
\ref{ss:SeparatedCovers},
we provide evidence of the fact that the set of non-metabelian groups \(G\)
whose second derived quotient \(G/G^{\prime\prime}\) is isomorphic to an assigned metabelian group \(M\)
is contained in the same tree as \(M\) itself, i.e.,
the \textit{covers} \(\mathrm{cov}(M)\) are separated by descendant trees.
In \S\
\ref{ss:Constraints}
we collect information on the relation rank of \(G\)
and the Galois action of the absolute group \(\mathrm{Gal}(K/\mathbb{Q})\) on \(G\),
which admits the determination of the length \(\ell_3(K)\)
of the \(3\)-class field tower of \(K\) in \S\
\ref{s:TowerLength}.


\section{Kernels of Artin transfers and abelian quotient invariants}
\label{s:ArtinPattern}


\subsection{Low index subgroups}
\label{ss:TopEndSbgrp}

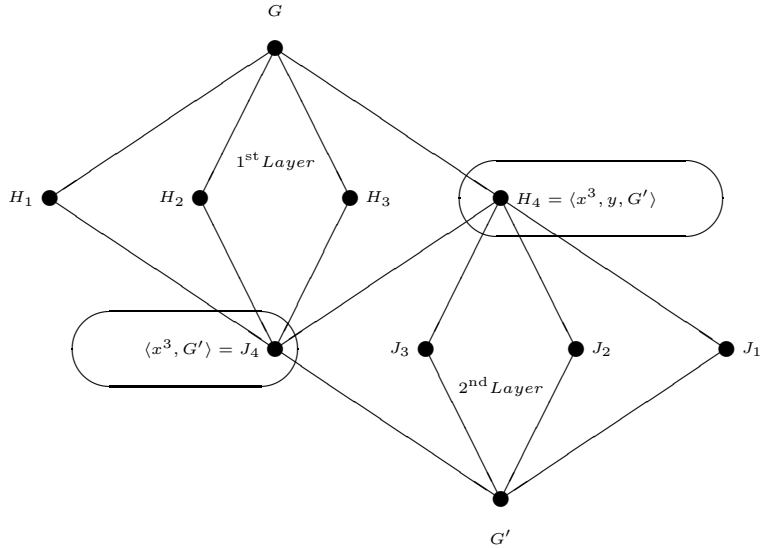
\begin{figure}[ht]
\caption{Low index subgroups of a pro-\(3\) group \(G\) with \(G/G^\prime\) of type \((9,3)\)}
\label{fig:DoubleDiamond}

{\tiny

\setlength{\unitlength}{1cm}
\begin{picture}(12,8)(-4.5,-7)

\put(0,0){\circle*{0.2}}
\put(0,0.5){\makebox(0,0)[cc]{\(G\)}}

\multiput(0,0)(3,-2){3}{\line(-3,-2){3}}
\multiput(0,0)(-3,-2){2}{\line(3,-2){3}}
\multiput(0,0)(1,-2){2}{\line(-1,-2){1}}

\put(0,-1.5){\makebox(0,0)[cc]{\(1^{\mathrm{st}} Layer\)}}
\multiput(-3,-2)(2,0){4}{\circle*{0.2}}
\put(-3.2,-2){\makebox(0,0)[rc]{\(H_1\)}}
\put(-1.2,-2){\makebox(0,0)[rc]{\(H_2\)}}
\put(1.2,-2){\makebox(0,0)[lc]{\(H_3\)}}
\put(3.2,-2){\makebox(0,0)[lc]{\(H_4=\langle x^3,y,G^\prime\rangle\)}}
\put(4.2,-2){\oval(3.5,1)}

\multiput(3,-2)(-3,-2){2}{\line(3,-2){3}}
\multiput(-1,-2)(1,2){2}{\line(1,-2){1}}
\multiput(3,-2)(1,-2){2}{\line(-1,-2){1}}

\put(3,-4.5){\makebox(0,0)[cc]{\(2^{\mathrm{nd}} Layer\)}}
\multiput(0,-4)(2,0){4}{\circle*{0.2}}
\put(6.2,-4){\makebox(0,0)[lc]{\(J_1\)}}
\put(4.2,-4){\makebox(0,0)[lc]{\(J_2\)}}
\put(1.8,-4){\makebox(0,0)[rc]{\(J_3\)}}
\put(-0.2,-4){\makebox(0,0)[rc]{\(\langle x^3,G^\prime\rangle=J_4\)}}
\put(-1.2,-4){\oval(3,1)}

\multiput(2,-4)(1,2){2}{\line(1,-2){1}}

\put(3,-6){\circle*{0.2}}
\put(3,-6.5){\makebox(0,0)[cc]{\(G^\prime\)}}

\end{picture}

}

\end{figure}

\noindent
Every two-generated pro-\(3\) group \(G=\langle x,y\rangle\)
with \(G/G^\prime\simeq (9,3)\),
such that \(x^9\in G^\prime\) and \(y^3\in G^\prime\),
possesses the following self-conjugate intermediate groups \(G^\prime\le J_i,H_i\le G\)
between the commutator subgroup \(G^\prime\) and \(G\),
as shown in Figure
\ref{fig:DoubleDiamond}:

\begin{itemize}

\item
\textit{first layer}: four \textit{maximal} normal subgroups \(H_i\) of index \((G:H_i)=3\),
\[
H_1=\langle x,G^\prime\rangle,\quad H_2=\langle xy,G^\prime\rangle,\quad
H_3=\langle xy^2,G^\prime\rangle,\quad H_4=\langle x^3,y,G^\prime\rangle, 
\]

\item
\textit{second layer}: four \textit{second maximal} normal subgroups \(J_i\) of index \((G:J_i)=9\),
\[
J_1=\langle y,G^\prime\rangle,\quad J_2=\langle x^3y,G^\prime\rangle,\quad
J_3=\langle x^3y^2,G^\prime\rangle,\quad J_4=\langle x^3,G^\prime\rangle.
\]

\end{itemize}


For both layers of subgroups, we use the subscript \(4\) to indicate
the \textit{distinguished} maximal subgroup, \(H_4=\prod_{i=1}^4\,J_i\),
for which the quotient \(H_4/G^\prime=\langle x^3,y\rangle\) is bicyclic of type \((3,3)\),
whereas \(H_i/G^\prime\) is cyclic of order \(9\) for \(1\le i\le 3\),
and to emphasize the \textit{distinguished} second maximal subgroup, \(J_4=\cap_{i=1}^4\,H_i\),
which coincides with the \textit{Frattini subgroup} \(\Phi(G)=G^3G^\prime\) of \(G\),
whereas \(J_i\) is contained in \(H_4\) only, for \(1\le i\le 3\).


\subsection{Artin transfers}
\label{ss:Transfers}

\noindent
We characterize a finite \(3\)-group \(G\)
by the usual invariants,
order \(\lvert G\rvert=3^m\),
logarithmic order \(\mathrm{lo}(G)=m\),
nilpotency class \(\mathrm{cl}(G)=c\),
coclass \(\mathrm{cc}(G)=m-c\),
and derived length \(\mathrm{dl}(G)=\ell\).
Additionally we use advanced invariants associated with
certain homomorphisms, the \textit{Artin transfers}
\cite{Ar1929}
from \(G\) to its maximal subgroups \(H_i\), \(1\le i\le 4\),
\[
T_i : G/G^\prime \longrightarrow H_i/H_i^\prime,\ g\mapsto
\begin{cases}
g^3                 & \text{ if } g\in G\setminus H_i \text{ (\textit{outer} transfer)}, \\
g^{\mathrm{S}_3(h)} & \text{ if } g\in H_i \text{ (\textit{inner} transfer)},
\end{cases}
\]
where \(\mathrm{S}_3(h)=1+h+h^2\in\mathbb{Z}\lbrack G\rbrack\),
with an arbitrary element \(h\in G\setminus H_i\),
denotes the third \textit{trace} element (\textit{Spur}) in the group ring of \(G\),
acting as symbolic exponent.
It turns out that, for fixed derived length \(\mathrm{dl}(G)=2\),
a metabelian group \(G\) is occasionally determined uniquely by
its \textit{transfer kernel type} (TKT), \(\varkappa(G)=(\ker(T_i))_{1\le i\le 4}\),
in conjunction with
its \textit{abelian quotient invariants} (AQI), \(\alpha(G)=(H_i/H_i^\prime)_{1\le i\le 4}\).
The \textit{Artin pattern} of \(G\) is the pair \(\mathrm{AP}(G)=(\alpha,\varkappa)\).
Here, we restrict the TKT and AQI to the \textit{first layer}.


\subsection{Orbits of punctured transfer kernel types}
\label{ss:Orbits}

\noindent
Although the capitulation over a few imaginary quadratic fields \(K\) of type \((9,3)\)
has been investigated in
\cite{SoTa1934,HeSm1982,Ma1991,Bb2012}
already, an \textit{invariant} characterization of the possible TKTs \(\varkappa(G_2)\) of their
second \(3\)-class groups \(G_2=\mathrm{Gal}(\mathrm{F}_3^2(K)/K)\) was missing up to now.
An adequate model, motivated by Figure
\ref{fig:DoubleDiamond},
is therefore established in the sequel, for the first time.

\begin{definition}
\label{dfn:TKT}
There are five possibilities for the kernel of the transfer \(T_i\), for each \(1\le i\le 4\).
Either \(\ker(T_i)=J_k/G^\prime\simeq C_3\), for some \(1\le k\le 4\),
and we denote the \textit{one-dimensional} transfer kernel by the singlet \(\varkappa(i)=k\)
\cite[\S\ 2.2, p. 475]{Ma2012a},
or \(\ker(T_i)=H_4/G^\prime\simeq C_3\times C_3\),
and we denote the \textit{two-dimensional} transfer kernel by the singlet \(\varkappa(i)=0\).
Due to the distinguished role of the subscript \(4\),
we combine the singlets in the following way to form a multiplet
\[\varkappa=\bigl(\,(\varkappa(1),\varkappa(2),\varkappa(3));\,\varkappa(4)\,\bigr)\in\lbrack 0,4\rbrack^3\times\lbrack 0,4\rbrack\]
which we call the \textit{punctured transfer kernel type} (pTKT) of the group \(G=\langle x,y\rangle\)
with respect to the selected generators \(x,y\).
In order to be independent of the choice of generators and of the arrangement of the subgroups
\(H_1,H_2,H_3\) and \(J_1,J_2,J_3\),
we define the \((S_3\times S_3)\)-\textit{orbit}
\[\varkappa^{S_3\times S_3}=\lbrace\tilde{\sigma}\circ\varkappa\circ\hat{\tau}\mid\sigma,\tau\in S_3\rbrace\]
of \(\varkappa\) under the operation of \(S_3\times S_3\)
as an \textit{isomorphism invariant} \(\varkappa(G)\) of \(G\).
Here, \(\tilde{\sigma}\) denotes the extension of \(\sigma\)
from \(\lbrack 1,3\rbrack\) to \(\lbrack 0,4\rbrack\)
which fixes \(0\) and \(4\),
and \(\hat{\tau}\) denotes the extension of \(\tau\)
from \(\lbrack 1,3\rbrack\) to \(\lbrack 1,4\rbrack\)
which fixes \(4\).
(The broader context of this definition is explained in
\cite{Ma2016b}.)
\end{definition}


\renewcommand{\arraystretch}{1.1}
\begin{table}[ht]
\caption{\(20\) \((S_3\times S_3)\)-orbits of \(\varkappa\in\lbrack 1,4\rbrack^4\) with \(\nu=0\)}
\label{tab:PrtPuncTrf}
\begin{center}
\begin{tabular}{|rr|c|ccccc|c|}
\hline
      &     & repres.       & occupation       &         & Taussky    & charact.                    & cardinality           & realizing                      \\
 Sec. & Nr. & of orbit      & numbers          &         & type       & property                    & of orbit              & \(3\)-group                    \\
      &     & \(\varkappa\) & \(o(\varkappa)\) & \(\mu\) & \(\kappa\) &                   & \(|\varkappa^{S_3\times S_3}|\) & \(G\)                          \\
\hline
    A &   1 & \((111;1)\)   & \((04000)\)      &   \(0\) & \((BBBA)\) & constant                    &   \(3\)               & \(\langle 3^4,6\rangle\)       \\
\hline
    B &   2 & \((111;2)\)   & \((03100)\)      &   \(0\) & \((BBBA)\) & nearly                      &   \(6\)               & \(\langle 3^8,1682\vert 1685\rangle\) \\
    B &   3 & \((112;1)\)   & \((03100)\)      &   \(0\) & \((BBBA)\) & constant                    &  \(18\)               &                                \\
\hline
    C &   4 & \((112;2)\)   & \((02200)\)      &   \(0\) & \((BBBA)\) &                             &  \(18\)               & \(\langle 3^8,1683\vert 1687\rangle\) \\
\hline
    D &   5 & \((112;3)\)   & \((02110)\)      &   \(0\) & \((BBBA)\) &                             &  \(18\)               & \(\langle 3^8,1684\vert 1686\rangle\) \\
    D &   6 & \((123;1)\)   & \((02110)\)      &   \(0\) & \((BBBA)\) &                             &  \(18\)               & \(\langle 3^8,1744\vert 1782\rangle\) \\
\hline
    B &   7 & \((111;4)\)   & \((03001)\)      &   \(1\) & \((BBBA)\) & nearly                      &   \(3\)               & \(\langle 3^6,16\vert 19\rangle\) \\
    B &   8 & \((114;1)\)   & \((03001)\)      &   \(1\) & \((BBAA)\) & constant                    &   \(9\)               &                                \\
\hline
    D &   9 & \((112;4)\)   & \((02101)\)      &   \(1\) & \((BBBA)\) &                             &  \(18\)               & exists                         \\
    D &  10 & \((114;2)\)   & \((02101)\)      &   \(1\) & \((BBAA)\) &                             &  \(18\)               & \(\langle 3^8,1689\vert 1690\rangle\) \\
    D &  11 & \((124;1)\)   & \((02101)\)      &   \(1\) & \((BBAA)\) &                             &  \(36\)               & \(\langle 3^6,14\vert 15\rangle\) \\
\hline
    E &  12 & \((123;4)\)   & \((01111)\)      &   \(1\) & \((BBBA)\) & per-                        &   \(6\)               & \(\langle 3^6,17\vert 20\rangle\) \\
    E &  13 & \((124;3)\)   & \((01111)\)      &   \(1\) & \((BBAA)\) & mutation                    &  \(18\)               &                                \\
\hline
    C &  14 & \((114;4)\)   & \((02002)\)      &   \(2\) & \((BBAA)\) &                             &   \(9\)               &                                \\
    C &  15 & \((144;1)\)   & \((02002)\)      &   \(2\) & \((BAAA)\) &                             &   \(9\)               &                                \\
\hline
    D &  16 & \((124;4)\)   & \((01102)\)      &   \(2\) & \((BBAA)\) &                             &  \(18\)               & exists                         \\
    D &  17 & \((144;2)\)   & \((01102)\)      &   \(2\) & \((BAAA)\) &                             &  \(18\)               & exists                         \\
\hline
    B &  18 & \((144;4)\)   & \((01003)\)      &   \(3\) & \((BAAA)\) & nearly                      &   \(9\)               & exists                         \\
    B &  19 & \((444;1)\)   & \((01003)\)      &   \(3\) & \((AAAA)\) & constant                    &   \(3\)               & exists                         \\
\hline
    A &  20 & \((444;4)\)   & \((00004)\)      &   \(4\) & \((AAAA)\) & constant                    &   \(1\)               & \(\langle 3^4,4\rangle,\langle 3^5,22\rangle,\langle 3^6,12\rangle\) \\
\hline
      &     &               &                  &         &            &               Total number: & \(256\)               &                                \\
\hline
\end{tabular}
\end{center}
\end{table}


\noindent
Throughout this work, we adhere to the convention that
the subscript \(4\) is \textit{distinguished}
and \textit{invariant} under any permutation of the other subscripts
(\(\tau\) for the domain and \(\sigma\) for the codomain).

\begin{definition}
\label{dfn:MuNu}
Two further \textit{isomorphism invariants} of \(G\) are defined by
the \textit{number of distinguished} transfer kernels
\(\mu=\mu(G)=\#\lbrace 1\le i\le 4\mid\varkappa(i)=4\rbrace\),
and the \textit{number of two-dimensional} transfer kernels
\(\nu=\nu(G)=\#\lbrace 1\le i\le 4\mid\varkappa(i)=0\rbrace\).
\end{definition}


\renewcommand{\arraystretch}{1.1}
\begin{table}[ht]
\caption{\(32\) \((S_3\times S_3)\)-orbits of \(\varkappa\in\lbrack 0,4\rbrack^4\setminus\lbrack 1,4\rbrack^4\) with \(1\le\nu\le 4\)}
\label{tab:TotPuncTrf}
\begin{center}
\begin{tabular}{|rr|c|ccccc|c|}
\hline
      &     & repres.       & occupation       &         & Taussky    & charact.                    & cardinality           & realizing                      \\
 Sec. & Nr. & of orbit      & numbers          &         & type       & property                    & of orbit              & \(3\)-group                    \\
      &     & \(\varkappa\) & \(o(\varkappa)\) & \(\mu\) & \(\kappa\) &                   & \(|\varkappa^{S_3\times S_3}|\) & \(G\)                          \\
\hline
    a &   1 & \((000;0)\)   & \((40000)\)      &   \(0\) & \((AAAA)\) & constant                    &   \(1\)               & \(\langle 3^4,3\rangle,\langle 3^5,15\vert 17\rangle\) \\
\hline
    b &   2 & \((000;1)\)   & \((31000)\)      &   \(0\) & \((AAAA)\) & nearly                      &   \(3\)               & \(\langle 3^5,16\rangle\)      \\
    b &   3 & \((001;0)\)   & \((31000)\)      &   \(0\) & \((AABA)\) & constant                    &   \(9\)               & \(\langle 3^5,19\vert 20\rangle\) \\
\hline
    c &   4 & \((001;1)\)   & \((22000)\)      &   \(0\) & \((AABA)\) &                             &   \(9\)               &                                \\
    c &   5 & \((011;0)\)   & \((22000)\)      &   \(0\) & \((ABBA)\) &                             &   \(9\)               &                                \\
\hline
    d &   6 & \((001;2)\)   & \((21100)\)      &   \(0\) & \((AABA)\) &                             &  \(18\)               &                                \\
    d &   7 & \((012;0)\)   & \((21100)\)      &   \(0\) & \((ABBA)\) &                             &  \(18\)               &                                \\
\hline
    b &   8 & \((011;1)\)   & \((13000)\)      &   \(0\) & \((ABBA)\) & nearly                      &   \(9\)               &                                \\
    b &   9 & \((111;0)\)   & \((13000)\)      &   \(0\) & \((BBBA)\) & constant                    &   \(3\)               &                                \\
\hline
    d &  10 & \((011;2)\)   & \((12100)\)      &   \(0\) & \((ABBA)\) &                             &  \(18\)               & \(\langle 3^6,13\rangle\)      \\
    d &  11 & \((012;1)\)   & \((12100)\)      &   \(0\) & \((ABBA)\) &                             &  \(36\)               &                                \\
    d &  12 & \((112;0)\)   & \((12100)\)      &   \(0\) & \((BBBA)\) &                             &  \(18\)               &                                \\
\hline
    e &  13 & \((012;3)\)   & \((11110)\)      &   \(0\) & \((ABBA)\) & per-                        &  \(18\)               &                                \\
    e &  14 & \((123;0)\)   & \((11110)\)      &   \(0\) & \((BBBA)\) & mutation                    &   \(6\)               & \(\langle 3^6,18\vert 21\rangle\) \\
\hline
    b &  15 & \((000;4)\)   & \((30001)\)      &   \(1\) & \((AAAA)\) & nearly                      &   \(1\)               & \(\langle 3^5,13\vert 14\rangle,\langle 3^6,9\rangle\) \\
    b &  16 & \((004;0)\)   & \((30001)\)      &   \(1\) & \((AAAA)\) & constant                    &   \(3\)               & \(\langle 3^5,18\rangle\)      \\
\hline
    d &  17 & \((001;4)\)   & \((21001)\)      &   \(1\) & \((AABA)\) &                             &   \(9\)               & exists                         \\
    d &  18 & \((004;1)\)   & \((21001)\)      &   \(1\) & \((AAAA)\) &                             &   \(9\)               &                                \\
    d &  19 & \((014;0)\)   & \((21001)\)      &   \(1\) & \((ABAA)\) &                             &  \(18\)               &                                \\
\hline
    d &  20 & \((011;4)\)   & \((12001)\)      &   \(1\) & \((ABBA)\) &                             &   \(9\)               &                                \\
    d &  21 & \((014;1)\)   & \((12001)\)      &   \(1\) & \((ABAA)\) &                             &  \(18\)               &                                \\
    d &  22 & \((114;0)\)   & \((12001)\)      &   \(1\) & \((BBAA)\) &                             &   \(9\)               &                                \\
\hline
    e &  23 & \((012;4)\)   & \((11101)\)      &   \(1\) & \((ABBA)\) & per-                        &  \(18\)               &                                \\
    e &  24 & \((014;2)\)   & \((11101)\)      &   \(1\) & \((ABAA)\) & muta-                       &  \(36\)               &                                \\
    e &  25 & \((124;0)\)   & \((11101)\)      &   \(1\) & \((BBAA)\) & tion                        &  \(18\)               &                                \\
\hline
    c &  26 & \((004;4)\)   & \((20002)\)      &   \(2\) & \((AAAA)\) &                             &   \(3\)               &                                \\
    c &  27 & \((044;0)\)   & \((20002)\)      &   \(2\) & \((AAAA)\) &                             &   \(3\)               & \(\langle 3^6,11\rangle\)      \\
\hline
    d &  28 & \((014;4)\)   & \((11002)\)      &   \(2\) & \((ABAA)\) &                             &  \(18\)               &                                \\
    d &  29 & \((044;1)\)   & \((11002)\)      &   \(2\) & \((AAAA)\) &                             &   \(9\)               &                                \\
    d &  30 & \((144;0)\)   & \((11002)\)      &   \(2\) & \((BAAA)\) &                             &   \(9\)               & exists                         \\
\hline
    b &  31 & \((044;4)\)   & \((10003)\)      &   \(3\) & \((AAAA)\) & nearly                      &   \(3\)               & \(\langle 3^6,10\rangle\)      \\
    b &  32 & \((444;0)\)   & \((10003)\)      &   \(3\) & \((AAAA)\) & constant                    &   \(1\)               &                                \\
\hline
      &     &               &                  &         & Total number: &             \(625-256=\) & \(369\)               &                                \\
\hline
\end{tabular}
\end{center}
\end{table}


\subsection{Combinatorially possible punctured transfer kernel types}
\label{ss:pTKT}

\noindent
In this section, we arrange all combinatorially possible
\((S_3\times S_3)\)-orbits of the \(5^4\) punctured quartets
\(\varkappa\in\lbrack 0,4\rbrack^3\times\lbrack 0,4\rbrack\)
by increasing invariant \(0\le\mu\le 4\) and cardinality of the image.
Table \ref{tab:PrtPuncTrf} shows the punctured quartets with invariant \(\nu=0\),
and Table \ref{tab:TotPuncTrf} those with invariant \(1\le\nu\le 4\),
as possible pTKTs \(\varkappa(G)\)
of \(3\)-groups \(G\) with \(G/G^\prime\) of type \((9,3)\),
respectively \textit{punctured capitulation types} \(\varkappa(K)\)
of number fields \(K\) with \(3\)-class group \(\mathrm{Cl}_3(K)\) of type \((9,3)\),
according to Artin's reciprocity law
\cite[\S\ 2.3, pp. 476--478]{Ma2012a}.
The orbits are divided into \textit{sections} (Sec), denoted by letters,
and identified by ordinal \textit{numbers} (Nr).
Each orbit contains a canonical \textit{representative}.


Table \ref{tab:PrtPuncTrf} gives
a coarse classification into sections \(\mathrm{A}\) to \(\mathrm{E}\),
an identification by ordinal numbers \(1\) to \(20\),
and a set theoretic characterization.
Table \ref{tab:TotPuncTrf} gives
a coarse classification into sections \(\mathrm{a}\) to \(\mathrm{e}\),
an identification by ordinal numbers \(1\) to \(32\),
and a set theoretic characterization.

We denote by \(o(\varkappa)=(\#\varkappa^{-1}\lbrace i\rbrace)_{0\le i\le 4}\)
the family of \textit{occupation numbers} of the selected orbit representative \(\varkappa\) and by
\(\kappa\) the quartet of \textit{Taussky's coarse capitulation types} \(A\) and \(B\) \cite{Ta1970} associated with \(\varkappa\)
(that is, \(\kappa(i)=A\) if the meet \(H_i\cap\ker(T_i)>1\) is non-trivial, and \(\kappa(i)=B\) otherwise).

If an orbit \(\varkappa^{S_3\times S_3}\) can be realized as pTKT \(\varkappa(G)\),
then a suitable \(3\)-group \(G\) is given
by its identifier in the SmallGroups library
\cite{BEO2005}.
In contrast to
\cite[Tbl. 6--7, pp. 492--493]{Ma2012a},
we are unable to mark an orbit as \lq\lq impossible\rq\rq\
when no realization as pTKT is known until now,
since currently, in contrast to \(G/G^\prime\simeq (3,3)\),
we do not have parametrized power-commutator presentations
of all metabelian \(3\)-groups \(G\) with \(G/G^\prime\simeq (9,3)\).


\renewcommand{\arraystretch}{1.1}
\begin{table}[ht]
\caption{Finite \(3\)-groups \(G\) with \(G/G^\prime\simeq C_9\times C_3\) and low order}
\label{tbl:LowOrder93}
\begin{center}
\begin{tabular}{|r|r||c|c|l||r|r||c|}
\hline
 ord     & id     & \(\alpha\)             & \(\varkappa\) & TKT               &   \(n\) & \(d_2\) & Action \\
\hline
  \(27\) &  \(2\) & \((2,2,2;11)\)         &   \((000;0)\) & \(\mathrm{a}.1\)  &   \(1\) &   \(3\) & \(\langle 12,4\rangle\) \\
\hline
  \(81\) &  \(3\) & \((21,21,21;111)\)     &   \((000;0)\) & \(\mathrm{a}.1\)  &   \(3\) &   \(4\) & \(\langle 12,4\rangle\) \\
  \(81\) &  \(4\) & \((21,21,21;21)\)      &   \((444;4)\) & \(\mathrm{A}.20\) &   \(2\) &   \(3\) & \(\langle 6,1\rangle\) \\
  \(81\) &  \(6\) & \((3,3,3;21)\)         &   \((111;1)\) & \(\mathrm{A}.1\)  &   \(0\) &   \(2\) & \(\langle 6,1\rangle\) \\
\hline
 \(243\) & \(13\) & \((21,21,21;1111)\)    &   \((000;4)\) & \(\mathrm{b}.15\) &   \(1\) &   \(4\) & \(\langle 12,4\rangle\) \\
 \(243\) & \(14\) & \((21,21,21;211)\)     &   \((000;4)\) & \(\mathrm{b}.15\) &   \(1\) &   \(4\) & \(\langle 12,4\rangle\) \\
 \(243\) & \(15\) & \((21,21,21;211)\)     &   \((000;0)\) & \(\mathrm{a}.1\)  &   \(1\) &   \(4\) & \(\langle 12,4\rangle\) \\
 \(243\) & \(16\) & \((21,21,21;211)\)     &   \((000;3)\) & \(\mathrm{b}.2\)  &   \(0\) &   \(3\) & \(\langle 6,2\rangle\) \\
 \(243\) & \(17\) & \((21,21,211;111)\)    &   \((000;0)\) & \(\mathrm{a}.1\)  &   \(1\) &   \(4\) & \(\langle 4,2\rangle\) \\
 \(243\) & \(18\) & \((21,21,211;111)\)    &   \((004;0)\) & \(\mathrm{b}.16\) &   \(0\) &   \(3\) & \(\langle 2,1\rangle\) \\
 \(243\) & \(19\) & \((21,21,31;111)\)     &   \((001;0)\) & \(\mathrm{b}.3\)  &   \(0\) &   \(3\) & \(\langle 4,2\rangle\) \\
 \(243\) & \(20\) & \((21,21,31;111)\)     &   \((001;0)\) & \(\mathrm{b}.3\)  &   \(0\) &   \(3\) & \(\langle 4,2\rangle\) \\
 \(243\) & \(22\) & \((21,21,21;21)\)      &   \((444;4)\) & \(\mathrm{A}.20\) &   \(0\) &   \(2\) & \(\langle 2,1\rangle\) \\
\hline
 \(729\) &  \(9\) & \((211,211,211;1111)\) &   \((000;4)\) & \(\mathrm{b}.15\) &   \(3\) &   \(5\) & \(\langle 12,4\rangle\) \\
 \(729\) & \(10\) & \((211,211,211;211)\)  &   \((044;4)\) & \(\mathrm{b}.31\) &   \(2\) &   \(4\) & \(\langle 4,2\rangle\) \\
 \(729\) & \(11\) & \((211,211,211;211)\)  &   \((044;0)\) & \(\mathrm{c}.27\) &   \(2\) &   \(4\) & \(\langle 4,2\rangle\) \\
 \(729\) & \(12\) & \((211,211,211;1111)\) &   \((444;4)\) & \(\mathrm{A}.20\) &   \(2\) &   \(4\) & \(\langle 6,2\rangle\) \\
\hline
 \(729\) & \(13\) & \((211,31,31;211)\)    &   \((011;3)\) & \(\mathrm{d}.10\) &   \(1\) &   \(3\) & \(\langle 2,1\rangle\) \\
 \(729\) & \(14\) & \((211,31,31;211)\)    &   \((423;3)\) & \(\mathrm{D}.11\) &   \(0\) &   \(2\) & \(\langle 2,1\rangle\) \\
 \(729\) & \(15\) & \((211,31,31;211)\)    &   \((432;3)\) & \(\mathrm{D}.11\) &   \(0\) &   \(2\) & \(\langle 2,1\rangle\) \\
 \(729\) & \(16\) & \((31,31,31;1111)\)    &   \((111;4)\) & \(\mathrm{B}.7\)  &   \(1\) &   \(3\) & \(\langle 12,4\rangle\) \\
 \(729\) & \(17\) & \((31,31,31;211)\)     &   \((123;4)\) & \(\mathrm{E}.12\) &   \(1\) &   \(3\) & \(\langle 12,4\rangle\) \\
 \(729\) & \(18\) & \((31,31,31;211)\)     &   \((132;0)\) & \(\mathrm{e}.14\) &   \(1\) &   \(3\) & \(\langle 12,4\rangle\) \\
 \(729\) & \(19\) & \((31,31,31;1111)\)    &   \((111;4)\) & \(\mathrm{B}.7\)  &   \(1\) &   \(3\) & \(\langle 12,4\rangle\) \\
 \(729\) & \(20\) & \((31,31,31;211)\)     &   \((132;4)\) & \(\mathrm{E}.12\) &   \(1\) &   \(3\) & \(\langle 12,4\rangle\) \\
 \(729\) & \(21\) & \((31,31,31;211)\)     &   \((123;0)\) & \(\mathrm{e}.14\) &   \(1\) &   \(3\) & \(\langle 12,4\rangle\) \\
\hline
\end{tabular}
\end{center}
\end{table}


\section{Finite \(3\)-groups with commutator quotient \((9,3)\)}
\label{s:NineByThree}

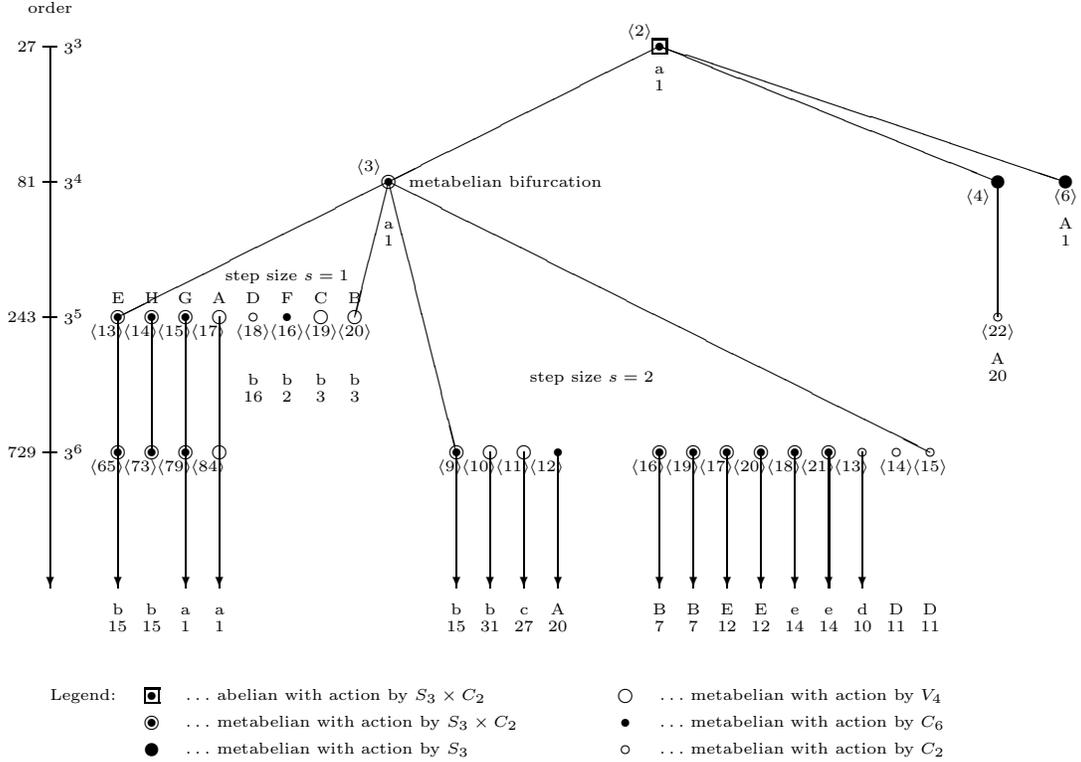
\begin{figure}[ht]
\caption{Finite \(3\)-groups \(G\) with commutator quotient \(G/G^\prime\simeq C_9\times C_3\)}
\label{fig:C9xC3}

{\tiny

\setlength{\unitlength}{0.9cm}
\begin{picture}(14,11.5)(-10,-10.5)

\put(-10,0.5){\makebox(0,0)[cb]{order}}

\put(-10,0){\line(0,-1){6}}
\multiput(-10.1,0)(0,-2){4}{\line(1,0){0.2}}

\put(-10.2,0){\makebox(0,0)[rc]{\(27\)}}
\put(-9.8,0){\makebox(0,0)[lc]{\(3^3\)}}
\put(-10.2,-2){\makebox(0,0)[rc]{\(81\)}}
\put(-9.8,-2){\makebox(0,0)[lc]{\(3^4\)}}
\put(-10.2,-4){\makebox(0,0)[rc]{\(243\)}}
\put(-9.8,-4){\makebox(0,0)[lc]{\(3^5\)}}
\put(-10.2,-6){\makebox(0,0)[rc]{\(729\)}}
\put(-9.8,-6){\makebox(0,0)[lc]{\(3^6\)}}

\put(-10,-6){\vector(0,-1){2}}

\put(-1.1,-0.1){\framebox(0.2,0.2){}}
\put(-1,0){\circle*{0.1}}


\multiput(-9,-4)(0.5,0){4}{\circle{0.2}}
\multiput(-9,-4)(0.5,0){3}{\circle*{0.1}}
\multiput(-9,-6)(0.5,0){4}{\circle{0.2}}
\multiput(-9,-6)(0.5,0){3}{\circle*{0.1}}

\put(-7,-4){\circle{0.1}}
\put(-6.5,-4){\circle*{0.1}}
\multiput(-6,-4)(0.5,0){2}{\circle{0.2}}

\put(-5,-2){\circle{0.2}}
\put(-5,-2){\circle*{0.1}}
\multiput(4,-2)(1,0){2}{\circle*{0.2}}
\put(4,-4){\circle{0.1}}

\put(-4,-6){\circle{0.2}}
\put(-4,-6){\circle*{0.1}}
\multiput(-3.5,-6)(0.5,0){2}{\circle{0.2}}
\put(-2.5,-6){\circle*{0.1}}

\multiput(-1,-6)(0.5,0){6}{\circle{0.2}}
\multiput(-1,-6)(0.5,0){6}{\circle*{0.1}}
\multiput(2,-6)(0.5,0){3}{\circle{0.1}}


\multiput(-9,-4)(0.5,0){4}{\line(0,-1){2}}
\put(-9,-6){\vector(0,-1){2}}
\multiput(-8,-6)(0.5,0){2}{\vector(0,-1){2}}

\put(-5,-2){\line(-2,-1){4}}
\put(-5,-2){\line(-1,-4){0.5}}
\put(-5,-2){\line(1,-4){1}}
\put(-5,-2){\line(2,-1){8}}

\put(-1,0){\line(-2,-1){4}}
\put(-1,0){\line(5,-2){5}}
\put(-1,0){\line(3,-1){6}}

\put(4,-2){\line(0,-1){2}}

\multiput(-4,-6)(0.5,0){4}{\vector(0,-1){2}}
\multiput(-1,-6)(0.5,0){7}{\vector(0,-1){2}}


\put(-1.1,0.1){\makebox(0,0)[rb]{\(\langle 2\rangle\)}}

\put(-6.5,-3.5){\makebox(0,0)[cb]{step size \(s=1\)}}
\put(-9,-3.8){\makebox(0,0)[cb]{\(\mathrm{E}\)}}
\put(-8.5,-3.8){\makebox(0,0)[cb]{\(\mathrm{H}\)}}
\put(-8,-3.8){\makebox(0,0)[cb]{\(\mathrm{G}\)}}
\put(-7.5,-3.8){\makebox(0,0)[cb]{\(\mathrm{A}\)}}
\put(-7,-3.8){\makebox(0,0)[cb]{\(\mathrm{D}\)}}
\put(-6.5,-3.8){\makebox(0,0)[cb]{\(\mathrm{F}\)}}
\put(-6,-3.8){\makebox(0,0)[cb]{\(\mathrm{C}\)}}
\put(-5.5,-3.8){\makebox(0,0)[cb]{\(\mathrm{B}\)}}

\put(-8.9,-4.1){\makebox(0,0)[rt]{\(\langle 13\rangle\)}}
\put(-8.9,-6.1){\makebox(0,0)[rt]{\(\langle 65\rangle\)}}
\put(-8.4,-4.1){\makebox(0,0)[rt]{\(\langle 14\rangle\)}}
\put(-8.4,-6.1){\makebox(0,0)[rt]{\(\langle 73\rangle\)}}
\put(-7.9,-4.1){\makebox(0,0)[rt]{\(\langle 15\rangle\)}}
\put(-7.9,-6.1){\makebox(0,0)[rt]{\(\langle 79\rangle\)}}
\put(-7.4,-4.1){\makebox(0,0)[rt]{\(\langle 17\rangle\)}}
\put(-7.4,-6.1){\makebox(0,0)[rt]{\(\langle 84\rangle\)}}
\put(-7,-4.1){\makebox(0,0)[ct]{\(\langle 18\rangle\)}}
\put(-6.5,-4.1){\makebox(0,0)[ct]{\(\langle 16\rangle\)}}
\put(-6,-4.1){\makebox(0,0)[ct]{\(\langle 19\rangle\)}}
\put(-5.5,-4.1){\makebox(0,0)[ct]{\(\langle 20\rangle\)}}

\put(-4.7,-2){\makebox(0,0)[lc]{metabelian bifurcation}}
\put(-5.1,-1.9){\makebox(0,0)[rb]{\(\langle 3\rangle\)}}
\put(3.9,-2.1){\makebox(0,0)[rt]{\(\langle 4\rangle\)}}
\put(4,-4.1){\makebox(0,0)[ct]{\(\langle 22\rangle\)}}
\put(4.5,-4.1){\makebox(0,0)[ct]{}}
\put(3.9,-6.1){\makebox(0,0)[rt]{}}
\put(4.5,-6.1){\makebox(0,0)[ct]{}}
\put(5,-2.1){\makebox(0,0)[ct]{\(\langle 6\rangle\)}}

\put(-2,-5){\makebox(0,0)[cb]{step size \(s=2\)}}
\put(-3.9,-6.1){\makebox(0,0)[rt]{\(\langle 9\rangle\)}}
\put(-3.4,-6.1){\makebox(0,0)[rt]{\(\langle 10\rangle\)}}
\put(-2.9,-6.1){\makebox(0,0)[rt]{\(\langle 11\rangle\)}}
\put(-2.4,-6.1){\makebox(0,0)[rt]{\(\langle 12\rangle\)}}
\put(-0.9,-6.1){\makebox(0,0)[rt]{\(\langle 16\rangle\)}}
\put(-0.4,-6.1){\makebox(0,0)[rt]{\(\langle 19\rangle\)}}
\put(0.1,-6.1){\makebox(0,0)[rt]{\(\langle 17\rangle\)}}
\put(0.6,-6.1){\makebox(0,0)[rt]{\(\langle 20\rangle\)}}
\put(1.1,-6.1){\makebox(0,0)[rt]{\(\langle 18\rangle\)}}
\put(1.6,-6.1){\makebox(0,0)[rt]{\(\langle 21\rangle\)}}
\put(2.1,-6.1){\makebox(0,0)[rt]{\(\langle 13\rangle\)}}
\put(2.5,-6.1){\makebox(0,0)[ct]{\(\langle 14\rangle\)}}
\put(3,-6.1){\makebox(0,0)[ct]{\(\langle 15\rangle\)}}


\put(-1,-0.4){\makebox(0,0)[cb]{\(\mathrm{a}\)}}
\put(-1,-0.5){\makebox(0,0)[ct]{\(1\)}}

\put(-9,-8.4){\makebox(0,0)[cb]{\(\mathrm{b}\)}}
\put(-9,-8.5){\makebox(0,0)[ct]{\(15\)}}
\put(-8.5,-8.4){\makebox(0,0)[cb]{\(\mathrm{b}\)}}
\put(-8.5,-8.5){\makebox(0,0)[ct]{\(15\)}}
\put(-8,-8.4){\makebox(0,0)[cb]{\(\mathrm{a}\)}}
\put(-8,-8.5){\makebox(0,0)[ct]{\(1\)}}
\put(-7.5,-8.4){\makebox(0,0)[cb]{\(\mathrm{a}\)}}
\put(-7.5,-8.5){\makebox(0,0)[ct]{\(1\)}}

\put(-7,-5){\makebox(0,0)[cb]{\(\mathrm{b}\)}}
\put(-7,-5.1){\makebox(0,0)[ct]{\(16\)}}
\put(-6.5,-5){\makebox(0,0)[cb]{\(\mathrm{b}\)}}
\put(-6.5,-5.1){\makebox(0,0)[ct]{\(2\)}}
\put(-6,-5){\makebox(0,0)[cb]{\(\mathrm{b}\)}}
\put(-6,-5.1){\makebox(0,0)[ct]{\(3\)}}
\put(-5.5,-5){\makebox(0,0)[cb]{\(\mathrm{b}\)}}
\put(-5.5,-5.1){\makebox(0,0)[ct]{\(3\)}}

\put(-5,-2.7){\makebox(0,0)[cb]{\(\mathrm{a}\)}}
\put(-5,-2.8){\makebox(0,0)[ct]{\(1\)}}

\put(4,-4.7){\makebox(0,0)[cb]{\(\mathrm{A}\)}}
\put(4,-4.8){\makebox(0,0)[ct]{\(20\)}}

\put(5,-2.7){\makebox(0,0)[cb]{\(\mathrm{A}\)}}
\put(5,-2.8){\makebox(0,0)[ct]{\(1\)}}

\put(-4,-8.4){\makebox(0,0)[cb]{\(\mathrm{b}\)}}
\put(-4,-8.5){\makebox(0,0)[ct]{\(15\)}}
\put(-3.5,-8.4){\makebox(0,0)[cb]{\(\mathrm{b}\)}}
\put(-3.5,-8.5){\makebox(0,0)[ct]{\(31\)}}
\put(-3,-8.4){\makebox(0,0)[cb]{\(\mathrm{c}\)}}
\put(-3,-8.5){\makebox(0,0)[ct]{\(27\)}}
\put(-2.5,-8.4){\makebox(0,0)[cb]{\(\mathrm{A}\)}}
\put(-2.5,-8.5){\makebox(0,0)[ct]{\(20\)}}

\put(-1,-8.4){\makebox(0,0)[cb]{\(\mathrm{B}\)}}
\put(-1,-8.5){\makebox(0,0)[ct]{\(7\)}}
\put(-0.5,-8.4){\makebox(0,0)[cb]{\(\mathrm{B}\)}}
\put(-0.5,-8.5){\makebox(0,0)[ct]{\(7\)}}
\put(0,-8.4){\makebox(0,0)[cb]{\(\mathrm{E}\)}}
\put(0,-8.5){\makebox(0,0)[ct]{\(12\)}}
\put(0.5,-8.4){\makebox(0,0)[cb]{\(\mathrm{E}\)}}
\put(0.5,-8.5){\makebox(0,0)[ct]{\(12\)}}
\put(1,-8.4){\makebox(0,0)[cb]{\(\mathrm{e}\)}}
\put(1,-8.5){\makebox(0,0)[ct]{\(14\)}}
\put(1.5,-8.4){\makebox(0,0)[cb]{\(\mathrm{e}\)}}
\put(1.5,-8.5){\makebox(0,0)[ct]{\(14\)}}
\put(2,-8.4){\makebox(0,0)[cb]{\(\mathrm{d}\)}}
\put(2,-8.5){\makebox(0,0)[ct]{\(10\)}}
\put(2.5,-8.4){\makebox(0,0)[cb]{\(\mathrm{D}\)}}
\put(2.5,-8.5){\makebox(0,0)[ct]{\(11\)}}
\put(3,-8.4){\makebox(0,0)[cb]{\(\mathrm{D}\)}}
\put(3,-8.5){\makebox(0,0)[ct]{\(11\)}}

\put(-10,-9.6){\makebox(0,0)[lc]{Legend:}}
\put(-8.6,-9.7){\framebox(0.2,0.2){}}
\put(-8.5,-9.6){\circle*{0.1}}
\put(-8,-9.6){\makebox(0,0)[lc]{\(\ldots\) abelian with action by \(S_3\times C_2\)}}
\put(-8.5,-10){\circle{0.2}}
\put(-8.5,-10){\circle*{0.1}}
\put(-8,-10){\makebox(0,0)[lc]{\(\ldots\) metabelian with action by \(S_3\times C_2\)}}
\put(-8.5,-10.4){\circle*{0.2}}
\put(-8,-10.4){\makebox(0,0)[lc]{\(\ldots\) metabelian with action by \(S_3\)}}
\put(-1.5,-9.6){\circle{0.2}}
\put(-1,-9.6){\makebox(0,0)[lc]{\(\ldots\) metabelian with action by \(V_4\)}}
\put(-1.5,-10){\circle*{0.1}}
\put(-1,-10){\makebox(0,0)[lc]{\(\ldots\) metabelian with action by \(C_6\)}}
\put(-1.5,-10.4){\circle{0.1}}
\put(-1,-10.4){\makebox(0,0)[lc]{\(\ldots\) metabelian with action by \(C_2\)}}

\end{picture}

}

\end{figure}

\noindent
In Table
\ref{tbl:LowOrder93}
we collect invariants of crucial metabelian \(3\)-groups \(G\) with \(G/G^\prime\simeq C_9\times C_3\)
and low order \(27\le\# G\le 729\).
The punctured AQI and TKT form the \textit{Artin pattern} \((\alpha,\varkappa)\) of \(G\),
which is used in the \textit{strategy of pattern recognition via Artin transfers}
\cite{Ma2020}
in order to identify the isomorphism class of the Galois group \(G_2=\mathrm{Gal}(\mathrm{F}_3^2(K)/K)\)
of the second Hilbert \(3\)-class field of a number field \(K\)
by the capitulation kernels \(\varkappa(K)\) and the abelian type invariants \(\alpha(K)\)
of the unramified cyclic cubic extensions \(L_i\), \(1\le i\le 4\),
of a number field with \(\mathrm{Cl}_3(K)\simeq C_9\times C_3\).
The \textit{nuclear rank} \(n\) specifies the position of \(G\) in a descendant tree, deciding
whether \(G\) is a terminal leaf with \(n=0\) or a root with further descendants if \(n\ge 1\).
The \textit{relation rank} \(d_2\) of \(G\) frequently admits an estimate of the
length of the Hilbert \(3\)-class field tower \(\mathrm{F}_3^\infty(K)\) of a number field \(K\).
Finally, the \textit{action} of the absolute Galois group \(\mathrm{Gal}(K/\mathbb{Q})\)
on the Frattini quotient \(Q=G/\Phi(G)\) decides
whether \(G\) is admissible as \(\mathrm{Gal}(\mathrm{F}_3^\infty(K)/K)\) for the field \(K\).
Generally, the isomorphism class of a group \(G\) is determined by its name in the
SmallGroups database
\cite{BEO2005}
which has the shape \(\langle\mathrm{ord},\mathrm{id}\rangle\)
containing the order and a numerical identifier in angle brackets.
Table
\ref{tbl:LowOrder93}
is illuminated by Figure
\ref{fig:C9xC3}.


By a \(3\)-group \textit{of type} \((9,3)\) we understand
a finite group \(G\) with derived quotient \(G/G^\prime\simeq C_9\times C_3\).
Such groups of second maximal class, that is of coclass \(\mathrm{cc}(G)=2\),
were called \textit{CF-groups} by Ascione et al.
\cite[\S\ 7, pp. 272--274]{AHL1977}.
Ascione denoted those of nilpotency class \(\mathrm{cl}(G)=3\)
by capital letters \(\mathrm{A},\ldots,\mathrm{H}\) as in Figure
\ref{fig:C9xC3}.
However, most of our \(3\)-class tower groups \(G=\mathrm{Gal}(\mathrm{F}_3^\infty(K)/K)\)
arise as descendants of \textit{step size} \(s=2\)
of the group \(\langle 81,3\rangle\)
with remarkable \textit{metabelian bifurcation} to coclass \(\mathrm{cc}=3\),
in the sense of
\cite{Ma2015a}.


\section{Finite \(3\)-groups with commutator quotient \((27,3)\) or \((81,3)\)}
\label{s:TwentysevenByThree}

\noindent
In Figure
\ref{fig:C27xC3},
resp.
\ref{fig:C81xC3},
we give identifiers of finite \(3\)-groups with commutator quotient
\(G/G^\prime\simeq (27,3)\), resp. \(G/G^\prime\simeq (81,3)\).
Directed edges lead from descendants \(D\) to parents \(\pi(D)=D/\gamma_c(D)\),
which usually differ from \(p\)-parents \(\pi_p(D)=D/P_{c_p-1}(D)\).
Up to minor modifications, the structure of the descendant trees in Figures
\ref{fig:C27xC3}
and
\ref{fig:C81xC3}
is the same as in Figure
\ref{fig:C9xC3}.
On the left hand side, there are the \textit{CF-groups} (cyclic factors) \(G\)
for which the factors \(\gamma_i(G)/\gamma_{i+1}(G)\simeq C_3\), \(i\ge 3\), are always cyclic.
They are higher analogues in branches \(\Phi_s(b)\), \(b\ge 1\),
of Ascione \(\mathrm{A},\ldots,\mathrm{H}\)
\cite{AHL1977}
in stems \(\Phi_s(0)\) of isoclinism classes.
On the right hand side, there are the \textit{non-CF groups} \(G\),
which are called \textit{BCF-groups} (bicyclic or cyclic factors) by Nebelung
\cite{Ne1989},
since \(\gamma_3(G)/\gamma_4(G)\simeq C_3\times C_3\) is bicyclic.


\begin{figure}[ht]
\caption{Finite \(3\)-groups \(G\) with commutator quotient \(G/G^\prime\simeq (27,3)\)}
\label{fig:C27xC3}



{\tiny

\setlength{\unitlength}{1.0cm}
\begin{picture}(14,9)(-9,-8)

\put(-10,0.5){\makebox(0,0)[cb]{order \(3^n\)}}

\put(-10,0){\line(0,-1){6}}
\multiput(-10.1,0)(0,-2){4}{\line(1,0){0.2}}

\put(-10.2,0){\makebox(0,0)[rc]{\(81\)}}
\put(-9.8,0){\makebox(0,0)[lc]{\(3^4\)}}
\put(-10.2,-2){\makebox(0,0)[rc]{\(243\)}}
\put(-9.8,-2){\makebox(0,0)[lc]{\(3^5\)}}
\put(-10.2,-4){\makebox(0,0)[rc]{\(729\)}}
\put(-9.8,-4){\makebox(0,0)[lc]{\(3^6\)}}
\put(-10.2,-6){\makebox(0,0)[rc]{\(2\,187\)}}
\put(-9.8,-6){\makebox(0,0)[lc]{\(3^7\)}}

\put(-10,-6){\vector(0,-1){2}}

\put(-8.6,0.2){\framebox(0.2,0.2){}}
\put(-8,0.3){\makebox(0,0)[lc]{\(\ldots\) abelian root}}
\put(-8.5,0){\circle{0.2}}
\put(-8,0){\makebox(0,0)[lc]{\(\ldots\) capable non-abelian}}
\put(-8.5,-0.3){\circle*{0.2}}
\put(-8,-0.3){\makebox(0,0)[lc]{\(\ldots\) terminal non-abelian}}

\put(-1.1,-0.1){\framebox(0.2,0.2){}}


\multiput(-9,-4)(0.5,0){4}{\circle{0.2}}
\multiput(-9,-6)(0.5,0){4}{\circle{0.2}}
\multiput(-7,-4)(0.5,0){4}{\circle*{0.2}}

\put(-5,-2){\circle{0.2}}
\put(4,-2){\circle{0.2}}
\put(4,-4){\circle{0.2}}
\put(4.5,-4){\circle*{0.2}}
\put(4,-6){\circle*{0.2}}
\put(5,-2){\circle*{0.2}}

\multiput(-4,-6)(0.5,0){4}{\circle{0.2}}
\multiput(-1,-6)(0.5,0){9}{\circle*{0.2}}


\multiput(-9,-4)(0.5,0){4}{\line(0,-1){2}}
\multiput(-9,-6)(0.5,0){4}{\vector(0,-1){2}}

\put(-5,-2){\line(-2,-1){4}}
\put(-5,-2){\line(-1,-4){0.5}}
\put(-5,-2){\line(1,-4){1}}
\put(-5,-2){\line(2,-1){8}}

\put(-1,0){\line(-2,-1){4}}
\put(-1,0){\line(5,-2){5}}
\put(-1,0){\line(3,-1){6}}

\put(4,-2){\line(0,-1){2}}
\put(4,-2){\line(1,-4){0.5}}
\put(4,-4){\line(0,-1){2}}

\multiput(-4,-6)(0.5,0){4}{\vector(0,-1){2}}


\put(-1.1,0.1){\makebox(0,0)[rb]{\(\langle 5\rangle\)}}

\put(-9,-4.1){\makebox(0,0)[rt]{\(\langle 4\rangle\)}}
\put(-8.5,-4.4){\makebox(0,0)[rt]{\(\langle 5\rangle\)}}
\put(-8,-4.1){\makebox(0,0)[rt]{\(\langle 6\rangle\)}}
\put(-7.5,-4.4){\makebox(0,0)[rt]{\(\langle 7\rangle\)}}
\put(-7,-4.1){\makebox(0,0)[ct]{\(\langle 8\rangle\)}}
\put(-6.5,-4.4){\makebox(0,0)[ct]{\(\langle 62\rangle\)}}
\put(-6,-4.1){\makebox(0,0)[ct]{\(\langle 63\rangle\)}}
\put(-5.5,-4.4){\makebox(0,0)[ct]{\(\langle 64\rangle\)}}

\put(-9,-6.1){\makebox(0,0)[rt]{\(\langle 86\rangle\)}}
\put(-9,-6.4){\makebox(0,0)[rt]{\(-\langle 92\rangle\)}}
\put(-8.5,-6.7){\makebox(0,0)[rt]{\(\langle 96\rangle\)}}
\put(-8.5,-7.0){\makebox(0,0)[rt]{\(-\langle 101\rangle\)}}
\put(-8,-6.1){\makebox(0,0)[rt]{\(\langle 105\rangle\)}}
\put(-8,-6.4){\makebox(0,0)[rt]{\(-\langle 110\rangle\)}}
\put(-7.5,-6.7){\makebox(0,0)[rt]{\(\langle 113\rangle\)}}
\put(-7.5,-7.0){\makebox(0,0)[rt]{\(-\langle 119\rangle\)}}

\put(-5.1,-1.9){\makebox(0,0)[rb]{\(\langle 12\rangle\)}}
\put(3.9,-2.1){\makebox(0,0)[rt]{\(\langle 21\rangle\)}}
\put(3.9,-4.1){\makebox(0,0)[rt]{\(\langle 22\rangle\)}}
\put(4.5,-4.1){\makebox(0,0)[ct]{\(\langle 92\rangle\)}}
\put(4,-6.1){\makebox(0,0)[ct]{\(\langle 194\rangle\)}}
\put(5,-2.1){\makebox(0,0)[ct]{\(\langle 24\rangle\)}}

\put(-4,-6.1){\makebox(0,0)[rt]{\(\langle 2\rangle\)}}
\put(-3.5,-6.4){\makebox(0,0)[rt]{\(\langle 3\rangle\)}}
\put(-3,-6.1){\makebox(0,0)[rt]{\(\langle 4\rangle\)}}
\put(-2.5,-6.4){\makebox(0,0)[rt]{\(\langle 5\rangle\)}}
\put(-1,-6.1){\makebox(0,0)[ct]{\(\langle 84\rangle\)}}
\put(-0.5,-6.4){\makebox(0,0)[ct]{\(\langle 85\rangle\)}}
\put(0,-6.1){\makebox(0,0)[ct]{\(\langle 94\rangle\)}}
\put(0.5,-6.4){\makebox(0,0)[ct]{\(\langle 95\rangle\)}}
\put(1,-6.1){\makebox(0,0)[ct]{\(\langle 103\rangle\)}}
\put(1.5,-6.4){\makebox(0,0)[ct]{\(\langle 104\rangle\)}}
\put(2,-6.1){\makebox(0,0)[ct]{\(\langle 112\rangle\)}}
\put(2.5,-6.4){\makebox(0,0)[ct]{\(\langle 121\rangle\)}}
\put(3,-6.1){\makebox(0,0)[ct]{\(\langle 122\rangle\)}}


\put(-1,-0.2){\makebox(0,0)[ct]{\(\mathrm{a}\)}}
\put(-1,-0.4){\makebox(0,0)[ct]{\(1\)}}

\put(-9,-8.1){\makebox(0,0)[ct]{\(\mathrm{b}\)}}
\put(-9,-8.3){\makebox(0,0)[ct]{\(15\)}}
\put(-8.5,-8.1){\makebox(0,0)[ct]{\(\mathrm{b}\)}}
\put(-8.5,-8.3){\makebox(0,0)[ct]{\(15\)}}
\put(-8,-8.1){\makebox(0,0)[ct]{\(\mathrm{a}\)}}
\put(-8,-8.3){\makebox(0,0)[ct]{\(1\)}}
\put(-7.5,-8.1){\makebox(0,0)[ct]{\(\mathrm{a}\)}}
\put(-7.5,-8.3){\makebox(0,0)[ct]{\(1\)}}

\put(-7,-4.7){\makebox(0,0)[ct]{\(\mathrm{b}\)}}
\put(-7,-4.9){\makebox(0,0)[ct]{\(16\)}}
\put(-6.5,-4.7){\makebox(0,0)[ct]{\(\mathrm{b}\)}}
\put(-6.5,-4.9){\makebox(0,0)[ct]{\(2\)}}
\put(-6,-4.7){\makebox(0,0)[ct]{\(\mathrm{b}\)}}
\put(-6,-4.9){\makebox(0,0)[ct]{\(3\)}}
\put(-5.5,-4.7){\makebox(0,0)[ct]{\(\mathrm{b}\)}}
\put(-5.5,-4.9){\makebox(0,0)[ct]{\(3\)}}

\put(-5,-2.3){\makebox(0,0)[ct]{\(\mathrm{a}\)}}
\put(-5,-2.5){\makebox(0,0)[ct]{\(1\)}}

\put(4,-6.3){\makebox(0,0)[ct]{\(\mathrm{A}\)}}
\put(4,-6.5){\makebox(0,0)[ct]{\(20\)}}

\put(5,-2.3){\makebox(0,0)[ct]{\(\mathrm{A}\)}}
\put(5,-2.5){\makebox(0,0)[ct]{\(1\)}}

\put(-4,-8.1){\makebox(0,0)[ct]{\(\mathrm{b}\)}}
\put(-4,-8.3){\makebox(0,0)[ct]{\(15\)}}
\put(-3.5,-8.1){\makebox(0,0)[ct]{\(\mathrm{b}\)}}
\put(-3.5,-8.3){\makebox(0,0)[ct]{\(31\)}}
\put(-3,-8.1){\makebox(0,0)[ct]{\(\mathrm{c}\)}}
\put(-3,-8.3){\makebox(0,0)[ct]{\(27\)}}
\put(-2.5,-8.1){\makebox(0,0)[ct]{\(\mathrm{A}\)}}
\put(-2.5,-8.3){\makebox(0,0)[ct]{\(20\)}}

\put(-1,-6.7){\makebox(0,0)[ct]{\(\mathrm{B}\)}}
\put(-1,-6.9){\makebox(0,0)[ct]{\(7\)}}
\put(-0.5,-6.7){\makebox(0,0)[ct]{\(\mathrm{B}\)}}
\put(-0.5,-6.9){\makebox(0,0)[ct]{\(7\)}}
\put(0,-6.7){\makebox(0,0)[ct]{\(\mathrm{E}\)}}
\put(0,-6.9){\makebox(0,0)[ct]{\(12\)}}
\put(0.5,-6.7){\makebox(0,0)[ct]{\(\mathrm{E}\)}}
\put(0.5,-6.9){\makebox(0,0)[ct]{\(12\)}}
\put(1,-6.7){\makebox(0,0)[ct]{\(\mathrm{e}\)}}
\put(1,-6.9){\makebox(0,0)[ct]{\(14\)}}
\put(1.5,-6.7){\makebox(0,0)[ct]{\(\mathrm{e}\)}}
\put(1.5,-6.9){\makebox(0,0)[ct]{\(14\)}}
\put(2,-6.7){\makebox(0,0)[ct]{\(\mathrm{d}\)}}
\put(2,-6.9){\makebox(0,0)[ct]{\(10\)}}
\put(2.5,-6.7){\makebox(0,0)[ct]{\(\mathrm{D}\)}}
\put(2.5,-6.9){\makebox(0,0)[ct]{\(11\)}}
\put(3,-6.7){\makebox(0,0)[ct]{\(\mathrm{D}\)}}
\put(3,-6.9){\makebox(0,0)[ct]{\(11\)}}

\end{picture}

}

\end{figure}


\begin{figure}[ht]
\caption{Finite \(3\)-groups \(G\) with commutator quotient \(G/G^\prime\simeq (81,3)\)}
\label{fig:C81xC3}



{\tiny

\setlength{\unitlength}{1.0cm}
\begin{picture}(14,9)(-9,-8)

\put(-10,0.5){\makebox(0,0)[cb]{order \(3^n\)}}

\put(-10,0){\line(0,-1){6}}
\multiput(-10.1,0)(0,-2){4}{\line(1,0){0.2}}

\put(-10.2,0){\makebox(0,0)[rc]{\(243\)}}
\put(-9.8,0){\makebox(0,0)[lc]{\(3^5\)}}
\put(-10.2,-2){\makebox(0,0)[rc]{\(729\)}}
\put(-9.8,-2){\makebox(0,0)[lc]{\(3^6\)}}
\put(-10.2,-4){\makebox(0,0)[rc]{\(2\,187\)}}
\put(-9.8,-4){\makebox(0,0)[lc]{\(3^7\)}}
\put(-10.2,-6){\makebox(0,0)[rc]{\(6\,561\)}}
\put(-9.8,-6){\makebox(0,0)[lc]{\(3^8\)}}

\put(-10,-6){\vector(0,-1){2}}

\put(-8.6,0.2){\framebox(0.2,0.2){}}
\put(-8,0.3){\makebox(0,0)[lc]{\(\ldots\) abelian root}}
\put(-8.5,0){\circle{0.2}}
\put(-8,0){\makebox(0,0)[lc]{\(\ldots\) capable non-abelian}}
\put(-8.5,-0.3){\circle*{0.2}}
\put(-8,-0.3){\makebox(0,0)[lc]{\(\ldots\) terminal non-abelian}}

\put(-1.1,-0.1){\framebox(0.2,0.2){}}


\multiput(-9,-4)(0.5,0){4}{\circle{0.2}}
\multiput(-9,-6)(0.5,0){4}{\circle{0.2}}
\multiput(-7,-4)(0.5,0){4}{\circle*{0.2}}

\put(-5,-2){\circle{0.2}}
\put(4,-2){\circle{0.2}}
\put(4,-4){\circle{0.2}}
\put(4.5,-4){\circle*{0.2}}
\put(4,-6){\circle{0.2}}
\put(4.5,-6){\circle*{0.2}}
\put(5,-2){\circle*{0.2}}

\multiput(-4,-6)(0.5,0){4}{\circle{0.2}}
\multiput(-1,-6)(0.5,0){9}{\circle*{0.2}}


\multiput(-9,-4)(0.5,0){4}{\line(0,-1){2}}
\multiput(-9,-6)(0.5,0){4}{\vector(0,-1){2}}

\put(-5,-2){\line(-2,-1){4}}
\put(-5,-2){\line(-1,-4){0.5}}
\put(-5,-2){\line(1,-4){1}}
\put(-5,-2){\line(2,-1){8}}

\put(-1,0){\line(-2,-1){4}}
\put(-1,0){\line(5,-2){5}}
\put(-1,0){\line(3,-1){6}}

\put(4,-2){\line(0,-1){2}}
\put(4,-2){\line(1,-4){0.5}}
\put(4,-4){\line(0,-1){2}}
\put(4,-4){\line(1,-4){0.5}}
\put(4,-6){\vector(0,-1){2}}

\multiput(-4,-6)(0.5,0){4}{\vector(0,-1){2}}


\put(-1.1,0.1){\makebox(0,0)[rb]{\(\langle 23\rangle\)}}

\put(-9,-4.1){\makebox(0,0)[rt]{\(\langle 83\rangle\)}}
\put(-8.5,-4.4){\makebox(0,0)[rt]{\(\langle 93\rangle\)}}
\put(-8,-4.1){\makebox(0,0)[rt]{\(\langle 102\rangle\)}}
\put(-7.5,-4.4){\makebox(0,0)[rt]{\(\langle 111\rangle\)}}
\put(-7,-4.1){\makebox(0,0)[ct]{\(\langle 120\rangle\)}}
\put(-6.5,-4.4){\makebox(0,0)[ct]{\(\langle 316\rangle\)}}
\put(-6,-4.1){\makebox(0,0)[ct]{\(\langle 317\rangle\)}}
\put(-5.5,-4.4){\makebox(0,0)[ct]{\(\langle 318\rangle\)}}

\put(-9,-6.1){\makebox(0,0)[rt]{\(\langle 63\rangle\)}}
\put(-9,-6.4){\makebox(0,0)[rt]{\(-\langle 66\rangle\)}}
\put(-8.5,-6.7){\makebox(0,0)[rt]{\(\langle 75\rangle\)}}
\put(-8.5,-7.0){\makebox(0,0)[rt]{\(-\langle 77\rangle\)}}
\put(-8,-6.1){\makebox(0,0)[rt]{\(\langle 84\rangle\)}}
\put(-8,-6.4){\makebox(0,0)[rt]{\(-\langle 86\rangle\)}}
\put(-7.5,-6.7){\makebox(0,0)[rt]{\(\langle 93\rangle\)}}
\put(-7.5,-7.0){\makebox(0,0)[rt]{\(-\langle 97\rangle\)}}

\put(-5.1,-1.9){\makebox(0,0)[rb]{\(\langle 61\rangle\)}}
\put(3.9,-2.1){\makebox(0,0)[rt]{\(\langle 91\rangle\)}}
\put(3.9,-4.1){\makebox(0,0)[rt]{\(\langle 193\rangle\)}}
\put(4.5,-4.1){\makebox(0,0)[ct]{\(\langle 383\rangle\)}}
\put(3.9,-6.1){\makebox(0,0)[rt]{\(\langle 199\rangle\)}}
\put(4.5,-6.1){\makebox(0,0)[ct]{\(\langle 1786\rangle\)}}
\put(5,-2.1){\makebox(0,0)[ct]{\(\langle 94\rangle\)}}

\put(-4,-6.1){\makebox(0,0)[rt]{\(\langle 200\rangle\)}}
\put(-3.5,-6.4){\makebox(0,0)[rt]{\(\langle 216\rangle\)}}
\put(-3,-6.1){\makebox(0,0)[rt]{\(\langle 229\rangle\)}}
\put(-2.5,-6.4){\makebox(0,0)[rt]{\(\langle 242\rangle\)}}
\put(-1,-6.1){\makebox(0,0)[ct]{\(\langle 876\rangle\)}}
\put(-0.5,-6.4){\makebox(0,0)[ct]{\(\langle 877\rangle\)}}
\put(0,-6.1){\makebox(0,0)[ct]{\(\langle 917\rangle\)}}
\put(0.5,-6.4){\makebox(0,0)[ct]{\(\langle 918\rangle\)}}
\put(1,-6.1){\makebox(0,0)[ct]{\(\langle 933\rangle\)}}
\put(1.5,-6.4){\makebox(0,0)[ct]{\(\langle 934\rangle\)}}
\put(2,-6.1){\makebox(0,0)[ct]{\(\langle 953\rangle\)}}
\put(2.5,-6.4){\makebox(0,0)[ct]{\(\langle 975\rangle\)}}
\put(3,-6.1){\makebox(0,0)[ct]{\(\langle 976\rangle\)}}


\put(-1,-0.2){\makebox(0,0)[ct]{\(\mathrm{a}\)}}
\put(-1,-0.4){\makebox(0,0)[ct]{\(1\)}}

\put(-9,-8.1){\makebox(0,0)[ct]{\(\mathrm{b}\)}}
\put(-9,-8.3){\makebox(0,0)[ct]{\(15\)}}
\put(-8.5,-8.1){\makebox(0,0)[ct]{\(\mathrm{b}\)}}
\put(-8.5,-8.3){\makebox(0,0)[ct]{\(15\)}}
\put(-8,-8.1){\makebox(0,0)[ct]{\(\mathrm{a}\)}}
\put(-8,-8.3){\makebox(0,0)[ct]{\(1\)}}
\put(-7.5,-8.1){\makebox(0,0)[ct]{\(\mathrm{a}\)}}
\put(-7.5,-8.3){\makebox(0,0)[ct]{\(1\)}}

\put(-7,-4.7){\makebox(0,0)[ct]{\(\mathrm{b}\)}}
\put(-7,-4.9){\makebox(0,0)[ct]{\(16\)}}
\put(-6.5,-4.7){\makebox(0,0)[ct]{\(\mathrm{b}\)}}
\put(-6.5,-4.9){\makebox(0,0)[ct]{\(2\)}}
\put(-6,-4.7){\makebox(0,0)[ct]{\(\mathrm{b}\)}}
\put(-6,-4.9){\makebox(0,0)[ct]{\(3\)}}
\put(-5.5,-4.7){\makebox(0,0)[ct]{\(\mathrm{b}\)}}
\put(-5.5,-4.9){\makebox(0,0)[ct]{\(3\)}}

\put(-5,-2.3){\makebox(0,0)[ct]{\(\mathrm{a}\)}}
\put(-5,-2.5){\makebox(0,0)[ct]{\(1\)}}

\put(4,-8.1){\makebox(0,0)[ct]{\(\mathrm{A}\)}}
\put(4,-8.3){\makebox(0,0)[ct]{\(20\)}}

\put(5,-2.3){\makebox(0,0)[ct]{\(\mathrm{A}\)}}
\put(5,-2.5){\makebox(0,0)[ct]{\(1\)}}

\put(-4,-8.1){\makebox(0,0)[ct]{\(\mathrm{b}\)}}
\put(-4,-8.3){\makebox(0,0)[ct]{\(15\)}}
\put(-3.5,-8.1){\makebox(0,0)[ct]{\(\mathrm{b}\)}}
\put(-3.5,-8.3){\makebox(0,0)[ct]{\(31\)}}
\put(-3,-8.1){\makebox(0,0)[ct]{\(\mathrm{c}\)}}
\put(-3,-8.3){\makebox(0,0)[ct]{\(27\)}}
\put(-2.5,-8.1){\makebox(0,0)[ct]{\(\mathrm{A}\)}}
\put(-2.5,-8.3){\makebox(0,0)[ct]{\(20\)}}

\put(-1,-6.7){\makebox(0,0)[ct]{\(\mathrm{B}\)}}
\put(-1,-6.9){\makebox(0,0)[ct]{\(7\)}}
\put(-0.5,-6.7){\makebox(0,0)[ct]{\(\mathrm{B}\)}}
\put(-0.5,-6.9){\makebox(0,0)[ct]{\(7\)}}
\put(0,-6.7){\makebox(0,0)[ct]{\(\mathrm{E}\)}}
\put(0,-6.9){\makebox(0,0)[ct]{\(12\)}}
\put(0.5,-6.7){\makebox(0,0)[ct]{\(\mathrm{E}\)}}
\put(0.5,-6.9){\makebox(0,0)[ct]{\(12\)}}
\put(1,-6.7){\makebox(0,0)[ct]{\(\mathrm{e}\)}}
\put(1,-6.9){\makebox(0,0)[ct]{\(14\)}}
\put(1.5,-6.7){\makebox(0,0)[ct]{\(\mathrm{e}\)}}
\put(1.5,-6.9){\makebox(0,0)[ct]{\(14\)}}
\put(2,-6.7){\makebox(0,0)[ct]{\(\mathrm{d}\)}}
\put(2,-6.9){\makebox(0,0)[ct]{\(10\)}}
\put(2.5,-6.7){\makebox(0,0)[ct]{\(\mathrm{D}\)}}
\put(2.5,-6.9){\makebox(0,0)[ct]{\(11\)}}
\put(3,-6.7){\makebox(0,0)[ct]{\(\mathrm{D}\)}}
\put(3,-6.9){\makebox(0,0)[ct]{\(11\)}}

\end{picture}

}

\end{figure}


\section{Length of the Hilbert \(3\)-class field tower}
\label{s:TowerLength}

\subsection{Separation of covers by descendant trees}
\label{ss:SeparatedCovers}

\noindent
In the following propositions,
let \(P\) and \(D\) be finite \(3\)-groups
with isomorphic commutator quotients \((9,3)\), and
let the parent \(P\) be a quotient of the descendant \(D\)
by a normal subgroup contained in the commutator subgroup \(D^\prime\).

\begin{proposition}
\label{prp:Quotients}
The components of \(\alpha(P)\) are quotients
of the corresponding components of \(\alpha(D)\),
which is exactly the meaning of the partial order relation \(\alpha(P)\le\alpha(D)\).
\end{proposition}

\begin{proof}
The statement is part of
the \textit{theorem on the antitony}
\(\alpha(P)\le\alpha(D)\) and \(\varkappa(P)\ge\varkappa(D)\)
of the components of the Artin pattern \((\alpha,\varkappa)\)
with respect to (parent, descendant)-pairs \((P,D)\),
where \(P\) is a quotient of \(D\),
which we proved in
\cite[\S\S\ 5.1--5.4, pp. 78--87]{Ma2016a}.
\end{proof}

\begin{corollary}
\label{cor:Quotients}
None of the metabelian descendants of
\(\langle 243,j\rangle\) with \(13\le j\le 20\)
can be the metabelianization of
any non-metabelian descendant of
\(\langle 729,i\rangle\) with \(13\le i\le 21\).
\end{corollary}

\begin{proof}
According to Table
\ref{tbl:LowOrder93},
the \(\alpha(N)\) of
non-metabelian descendants \(N\) of
the groups \(\langle 729,i\rangle\) with \(13\le i\le 21\)
have at least two components \((31)\).
Since the \(\alpha(M)\) of
metabelian descendants \(M\) of
the groups \(\langle 243,j\rangle\) with \(13\le j\le 20\)
have at least two components \((21)\),
none of the metabelianizations \(N/N^{\prime\prime}\)
(each of which has AQI coinciding with \(\alpha(N)\))
can be isomorphic to one of the groups \(M\).
\end{proof}


\begin{proposition}
\label{prp:Ranks}
The ranks of the components of \(\alpha(D)\) cannot be smaller than
the ranks of the corresponding components of \(\alpha(P)\).
\end{proposition}

\begin{proof}
This is an immediate consequence of Proposition
\ref{prp:Quotients}.
\end{proof}

\begin{corollary}
\label{cor:Ranks}
None of the metabelian descendants of
\(\langle 243,j\rangle\) with \(13\le j\le 20\)
can be the metabelianization of
any non-metabelian descendant of
\(\langle 729,i\rangle\) with \(9\le i\le 12\).
\end{corollary}

\begin{proof}
According to Table
\ref{tbl:LowOrder93},
the four components of \(\alpha(P)\) have at least rank three,
for the groups \(P=\langle 729,i\rangle\) with \(9\le i\le 12\).
By the \textit{antitony principle},
this is also true for \(\alpha(D)\)
of any non-metabelian descendant \(D\)
of one of these four roots \(P\).
However,
metabelian descendants \(M\) of
the groups \(\langle 243,j\rangle\) with \(13\le j\le 20\)
have \(\alpha(M)\) with at least two components of rank two.
Since each of the metabelianizations \(D/D^{\prime\prime}\)
has AQI coinciding with \(\alpha(D)\),
none of them can be isomorphic to one of the groups \(M\).
\end{proof}


\subsection{Relation rank and Galois action}
\label{ss:Constraints}

\noindent
Constraints arise from two issues,
bounds for the relation rank of the tower group \(G=\mathrm{Gal}(\mathrm{F}_3^\infty(K)/K)\),
and the Galois action of \(\mathrm{Gal}(K/\mathbb{Q})\) on \(\mathrm{Cl}_3(K)\simeq G/G^\prime\).
By \(\langle o,i\rangle\) we denote groups in the SmallGroups database of Magma
\cite{MAGMA2021}.
In tree diagrams, the order \(o\) is given on a scale, and we abbreviate the identifiers by \(\langle i\rangle\).


\begin{theorem}
\label{thm:RelationRankAndGaloisAction}
For a number field \(K\)
with \(3\)-class rank \(\varrho=2\),
in particular for \(\mathrm{Cl}_3(K)\simeq C_9\times C_3\),
the Galois group \(G=\mathrm{Gal}(\mathrm{F}_3^\infty(K)/K)\) of the \(3\)-class field tower
must satisfy the following conditions.
\begin{enumerate}
\item
The relation rank \(d_2\) of \(G\) must be bounded by
\(2\le d_2\le 2+r+\theta\),
where \(r=r_1+r_2-1\) denotes the torsion free Dirichlet unit rank
of the field \(K\) with signature \((r_1,r_2)\),
and \(\theta=1\), if \(K\) contains the primitive third roots of unity,
\(\theta=0\) otherwise.
\item
The automorphism group \(\mathrm{Aut}(Q)\) of the Frattini quotient \(Q=G/\Phi(G)\)
must contain a subgroup isomorphic to
\(\mathrm{Gal}(K/\mathbb{Q})\).
\end{enumerate}
\end{theorem}

\begin{proof}
According to the Burnside basis theorem,
the generator rank \(d_1\) of \(G\) coincides with
the generator rank of the Frattini quotient \(Q=G/\Phi(G)=G/(G^\prime\cdot G^3)\),
respectively the derived quotient \(G/G^\prime\simeq\mathrm{Cl}_3(K)\),
that is the \(3\)-class rank \(\varrho=2\) of \(K\).
\begin{enumerate}
\item
According to the Shafarevich Theorem
\cite[Thm. 5.1, p. 28]{Ma2015b},
the relation rank \(d_2\) of \(G\) is bounded by
\(d_1\le d_2\le d_1+r+\theta\).
Together with the generator rank \(d_1=\varrho=2\) this gives the bounds
\(2\le d_2\le 2+r+\theta\).
\item
The absolute Galois group \(\mathrm{Gal}(K/\mathbb{Q})\) of \(K\)
acts on the \(3\)-class group \(\mathrm{Cl}_3(K)\simeq G/G^\prime\)
and thus also on the Frattini quotient \(Q=G/\Phi(G)=G/(G^\prime\cdot G^3)\),
whence \(\mathrm{Aut}(Q)\) contains a subgroup isomorphic to \(\mathrm{Gal}(K/\mathbb{Q})\). \qedhere
\end{enumerate}
\end{proof}


\noindent
By the same proof as for item (2) of Theorem
\ref{thm:RelationRankAndGaloisAction},
with \(G/G^\prime\simeq\mathrm{Cl}_3(K)\)
replaced by
\[
G_k/G_k^\prime\simeq
\mathrm{Gal}(\mathrm{F}_3^k(K)/K)\bigm/\mathrm{Gal}(\mathrm{F}_3^k(K)/\mathrm{F}_3^1(K))\simeq
\mathrm{Gal}(\mathrm{F}_3^1(K)/K)\simeq\mathrm{Cl}_3(K)
\]
we obtain the same requirement for the Galois action on \(G_k\) (but \textit{not} for the relation rank of \(G_k\)):

\begin{corollary}
\label{cor:GaloisAction}
Let \(k\) be a positive integer,
and denote by \(G_k=\mathrm{Gal}(\mathrm{F}_3^k(K)/K)\)
the Galois group of the \(k\)-th Hilbert \(3\)-class field \(\mathrm{F}_3^k(K)\) of \(K\).
The automorphism group \(\mathrm{Aut}(Q)\) of the Frattini quotient \(Q=G_k/\Phi(G_k)\)
must contain a subgroup isomorphic to
\(\mathrm{Gal}(K/\mathbb{Q})\).
\end{corollary}


\begin{figure}[ht]
\caption{Metabelian skeleton of the coclass tree \(\mathcal{T}_3\langle 729,13\rangle\)}
\label{fig:Ord729Id13}

{\tiny

\setlength{\unitlength}{1cm}
\begin{picture}(14,14.5)(-6,-13.5)

\put(-5,0.5){\makebox(0,0)[cb]{order}}

\put(-5,0){\line(0,-1){12}}
\multiput(-5.1,0)(0,-2){7}{\line(1,0){0.2}}

\put(-5.2,0){\makebox(0,0)[rc]{\(729\)}}
\put(-4.8,0){\makebox(0,0)[lc]{\(3^6\)}}
\put(-5.2,-2){\makebox(0,0)[rc]{\(2\,187\)}}
\put(-4.8,-2){\makebox(0,0)[lc]{\(3^7\)}}
\put(-5.2,-4){\makebox(0,0)[rc]{\(6\,561\)}}
\put(-4.8,-4){\makebox(0,0)[lc]{\(3^8\)}}
\put(-5.2,-6){\makebox(0,0)[rc]{\(19\,683\)}}
\put(-4.8,-6){\makebox(0,0)[lc]{\(3^9\)}}
\put(-5.2,-8){\makebox(0,0)[rc]{\(59\,049\)}}
\put(-4.8,-8){\makebox(0,0)[lc]{\(3^{10}\)}}
\put(-5.2,-10){\makebox(0,0)[rc]{\(177\,147\)}}
\put(-4.8,-10){\makebox(0,0)[lc]{\(3^{11}\)}}
\put(-5.2,-12){\makebox(0,0)[rc]{\(531\,441\)}}
\put(-4.8,-12){\makebox(0,0)[lc]{\(3^{12}\)}}

\put(-5,-12){\vector(0,-1){2}}

\put(-3.5,0.3){\circle{0.2}}
\put(-3.5,0.3){\circle{0.1}}
\put(-3,0.3){\makebox(0,0)[lc]{\(\ldots\) GI action by \(C_2\), homocyclic AQI}}
\put(-3.5,0){\circle{0.2}}
\put(-3,0){\makebox(0,0)[lc]{\(\ldots\) GI and RI action by \(C_2\), homocyclic AQI}}

\put(4.5,0.3){\circle{0.2}}
\put(4.5,0.3){\circle*{0.1}}
\put(5,0.3){\makebox(0,0)[lc]{\(\ldots\) GI action by \(C_2\)}}
\put(4.5,0){\circle*{0.2}}
\put(5,0){\makebox(0,0)[lc]{\(\ldots\) GI and RI action by \(C_2\)}}
\put(4.5,-0.3){\circle{0.1}}
\put(5,-0.3){\makebox(0,0)[lc]{\(\ldots\) non-\(\sigma\)}}

\put(3.2,0){\makebox(0,0)[lb]{\(\langle 13\rangle\)}}
\put(3.2,-2){\makebox(0,0)[lb]{\(\langle 168\rangle\)}}
\put(3.1,-4){\makebox(0,0)[lb]{\(\langle 1688\rangle\)}}
\multiput(3.2,-6)(0,-2){4}{\makebox(0,0)[lc]{\(4\)}}

\multiput(3,0)(0,-4){4}{\circle{0.2}}
\multiput(3,0)(0,-4){4}{\circle*{0.1}}

\multiput(3,-2)(0,-4){3}{\circle*{0.2}}

\multiput(3,0)(0,-2){6}{\line(0,-1){2}}

\put(3,-12){\vector(0,-1){2}}
\put(3.2,-12.5){\makebox(0,0)[lc]{main line}}

\put(-2.8,-2){\makebox(0,0)[lc]{\(\langle 171\rangle\)}}
\put(-0.8,-2){\makebox(0,0)[lc]{\(\langle 172\rangle\)}}
\put(1.2,-2){\makebox(0,0)[lc]{\(\langle 169\rangle\)}}
\put(5.2,-2){\makebox(0,0)[lc]{\(\langle 170\rangle\)}}

\put(-2.9,-4){\makebox(0,0)[lc]{\(\langle 1683\rangle\)}}
\put(-1.9,-4){\makebox(0,0)[lc]{\(\langle 1687\rangle\)}}
\put(-0.9,-4){\makebox(0,0)[lc]{\(\langle 1684\rangle\)}}
\put(0.1,-4){\makebox(0,0)[lc]{\(\langle 1686\rangle\)}}
\put(1.1,-4){\makebox(0,0)[lc]{\(\langle 1689\rangle\)}}
\put(2.1,-4){\makebox(0,0)[lc]{\(\langle 1690\rangle\)}}
\put(5.1,-4){\makebox(0,0)[lc]{\(\langle 1682\rangle\)}}
\put(6.1,-4){\makebox(0,0)[lc]{\(\langle 1685\rangle\)}}

\multiput(-2.8,-6)(0,-4){2}{\makebox(0,0)[lc]{\(2\)}}
\multiput(-2.8,-8)(0,-4){2}{\makebox(0,0)[lc]{\(3\)}}
\multiput(-1.8,-8)(0,-4){2}{\makebox(0,0)[lc]{\(8\)}}
\multiput(-0.8,-6)(0,-4){2}{\makebox(0,0)[lc]{\(3\)}}
\multiput(-0.8,-8)(0,-4){2}{\makebox(0,0)[lc]{\(2\)}}
\multiput(0.2,-8)(0,-4){2}{\makebox(0,0)[lc]{\(9\)}}
\multiput(1.2,-6)(0,-2){4}{\makebox(0,0)[lc]{\(5\)}}
\multiput(2.2,-8)(0,-4){2}{\makebox(0,0)[lc]{\(6\)}}
\multiput(5.2,-6)(0,-2){4}{\makebox(0,0)[lc]{\(1\)}}
\multiput(6.2,-8)(0,-4){2}{\makebox(0,0)[lc]{\(7\)}}

\multiput(-3,-2)(0,-4){3}{\circle{0.1}}

\put(-3,-4){\circle{0.2}}
\multiput(-3,-8)(0,-4){2}{\circle*{0.2}}
\multiput(-2,-4)(0,-4){3}{\circle*{0.2}}

\multiput(-1,-2)(0,-4){3}{\circle{0.1}}

\put(-1,-4){\circle{0.2}}
\multiput(-1,-8)(0,-4){2}{\circle*{0.2}}
\multiput(0,-4)(0,-4){3}{\circle*{0.2}}

\multiput(1,-2)(0,-4){3}{\circle{0.1}}

\multiput(1,-4)(0,-4){3}{\circle*{0.2}}

\multiput(2,-4)(0,-4){3}{\circle*{0.2}}

\multiput(5,-2)(0,-4){3}{\circle{0.1}}

\multiput(5,-4)(0,-4){3}{\circle{0.2}}
\put(5,-4){\circle{0.1}}
\multiput(5,-8)(0,-4){2}{\circle*{0.1}}
\multiput(6,-4)(0,-4){3}{\circle{0.2}}
\multiput(6,-4)(0,-4){3}{\circle*{0.1}}

\multiput(3,0)(0,-2){6}{\line(-3,-1){6}}
\multiput(3,0)(0,-2){6}{\line(-2,-1){4}}
\multiput(3,0)(0,-2){6}{\line(-1,-1){2}}
\multiput(3,0)(0,-2){6}{\line(1,-1){2}}
\multiput(3,-2)(0,-4){3}{\line(3,-2){3}}

\multiput(5,-4)(0,-4){2}{\line(1,-1){2}}
\multiput(6,-4)(0,-4){2}{\line(1,-1){2}}

\put(-4,-13){\makebox(0,0)[cc]{pTKT}}
\put(-4,-13.5){\makebox(0,0)[cc]{\(\varkappa\sim\)}}
\put(-2.5,-13){\makebox(0,0)[cc]{\(\mathrm{C}.4\)}}
\put(-2.5,-13.5){\makebox(0,0)[cc]{\((311;3)\)}}
\put(-0.5,-13){\makebox(0,0)[cc]{\(\mathrm{D}.5\)}}
\put(-0.5,-13.5){\makebox(0,0)[cc]{\((211;3)\)}}
\put(1.5,-13){\makebox(0,0)[cc]{\(\mathrm{D}.10\)}}
\put(1.5,-13.5){\makebox(0,0)[cc]{\((411;3)\)}}
\put(3.5,-13){\makebox(0,0)[cc]{\(\mathrm{d}.10\)}}
\put(3.5,-13.5){\makebox(0,0)[cc]{\((011;3)\)}}
\put(5.5,-13){\makebox(0,0)[cc]{\(\mathrm{B}.2\)}}
\put(5.5,-13.5){\makebox(0,0)[cc]{\((111;3)\)}}

\end{picture}

}

\end{figure}
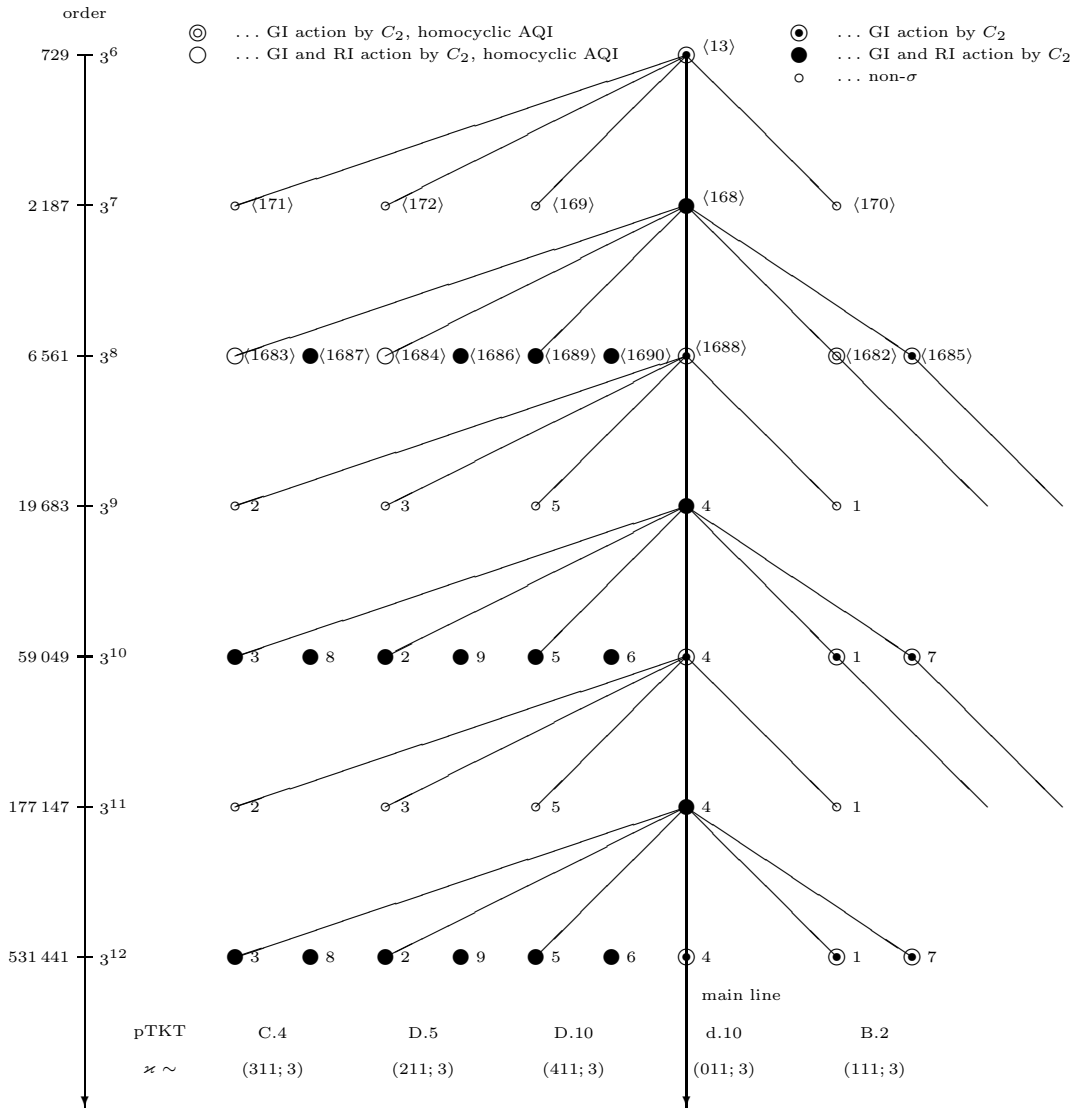


\section{Imaginary quadratic fields with \(3\)-class group \((9,3)\)}
\label{s:ImagQuadratic}

\noindent
An imaginary quadratic field
\(K=\mathbb{Q}(\sqrt{d})\)
has signature \((r_1,r_2)=(0,1)\),
and torsionfree Dirichlet unit rank \(r=r_1+r_2-1=0\).
If it has class number bigger than one,
then it does not contain
primitive third roots of unity, i.e. \(\theta=0\).
Shafarevich bounds for the relation rank of
\(G_\infty=\mathrm{Gal}(\mathrm{F}_3^\infty(K)/K)\) are given by
\(2=d_1\le d_2\le d_1+r+\theta=2\),
i.e., \textit{Schur \(\sigma\)-groups} are mandatory.
Among the \(875\) imaginary quadratic fields \(K=\mathbb{Q}(\sqrt{d})\)
with fundamental discriminants \(-10^6<d<0\)
and \(\mathrm{Cl}_3(K)\simeq C_9\times C_3\),
the pTKTs and second \(3\)-class groups are distributed as in Table
\ref{tbl:ImagQuadratic}.


\renewcommand{\arraystretch}{1.0}
\begin{table}[ht]
\caption{Length \(\ell_3(K)\) of \(3\)-class field towers of imaginary quadratic fields \(K\)}
\label{tbl:ImagQuadratic}
\begin{center}
\begin{tabular}{|r|r|r||l|c|c|c|}
\hline
 \(\#\)  & \(\%\)   & \(\lvert d_0\rvert\) & pTKT              & AQI          &  \(\mathrm{Gal}(\mathrm{F}_3^2(K)/K)\)            & \(\ell_3(K)\) \\
\hline
 \(406\) & \(46.4\) &           \(3\,299\) & \(\mathrm{D}.11\) &              &  \(\langle 729,14\vert 15\rangle\)                & \(2\) \\
\hline
  \(75\) &  \(8.6\) &           \(5\,703\) & \(\mathrm{E}.12\) &              &  \(\in\mathcal{T}\langle 729,17\vert 20\rangle\)  & \(\ge 3\) \\
  \(64\) &  \(7.3\) &          \(54\,695\) & \(\mathrm{B}.7\)  &              &  \(\in\mathcal{T}\langle 729,16\vert 19\rangle\)  & \(\ge 3\) \\
  \(20\) &  \(2.3\) &         \(289\,704\) & \(\mathrm{A}.20\) &              &  \(\in\mathcal{T}\langle 729,12\rangle\)          & \(\ge 3\) \\
\hline
  \(46\) &  \(5.3\) &          \(11\,651\) & \(\mathrm{D}.10\) &              &  \(\langle 6561,1689\vert 1690\rangle\)           & \(3\) \\
  \(20\) &  \(2.3\) &          \(17\,723\) & \(\mathrm{D}.5\)  & heterocyclic &  \(\langle 6561,1686\rangle\)                     & \(3\) \\
  \(40\) &  \(4.6\) &          \(31\,983\) & \(\mathrm{D}.6\)  &              &  \(\langle 6561,1744\vert 1782\rangle\)           & \(3\) \\
  \(14\) &  \(1.6\) &          \(35\,331\) & \(\mathrm{C}.4\)  & heterocyclic &  \(\langle 6561,1687\rangle\)                     & \(3\) \\
  \(30\) &  \(3.4\) &          \(42\,567\) & \(\mathrm{C}.4\)  & homocyclic   &  \(\langle 6561,1683\rangle\)                     & \(3\) \\
  \(15\) &  \(1.7\) &         \(116\,419\) & \(\mathrm{D}.5\)  & homocyclic   &  \(\langle 6561,1684\rangle\)                     & \(3\) \\
\hline
\end{tabular}
\end{center}
\end{table}


\begin{figure}[hb]
\caption{Metabelian skeleton of the coclass tree \(\mathcal{T}_3\langle 729,21\rangle\)}
\label{fig:Ord729Id21}

{\tiny

\setlength{\unitlength}{0.9cm}
\begin{picture}(14,14.8)(-4,-13.5)

\put(-2,0.5){\makebox(0,0)[cb]{order}}

\put(-2,0){\line(0,-1){12}}
\multiput(-2.1,0)(0,-2){7}{\line(1,0){0.2}}

\put(-2.2,0){\makebox(0,0)[rc]{\(729\)}}
\put(-1.8,0){\makebox(0,0)[lc]{\(3^6\)}}
\put(-2.2,-2){\makebox(0,0)[rc]{\(2\,187\)}}
\put(-1.8,-2){\makebox(0,0)[lc]{\(3^7\)}}
\put(-2.2,-4){\makebox(0,0)[rc]{\(6\,561\)}}
\put(-1.8,-4){\makebox(0,0)[lc]{\(3^8\)}}
\put(-2.2,-6){\makebox(0,0)[rc]{\(19\,683\)}}
\put(-1.8,-6){\makebox(0,0)[lc]{\(3^9\)}}
\put(-2.2,-8){\makebox(0,0)[rc]{\(59\,049\)}}
\put(-1.8,-8){\makebox(0,0)[lc]{\(3^{10}\)}}
\put(-2.2,-10){\makebox(0,0)[rc]{\(177\,147\)}}
\put(-1.8,-10){\makebox(0,0)[lc]{\(3^{11}\)}}
\put(-2.2,-12){\makebox(0,0)[rc]{\(531\,441\)}}
\put(-1.8,-12){\makebox(0,0)[lc]{\(3^{12}\)}}

\put(-2,-12){\vector(0,-1){2}}

\put(4.5,0.3){\circle{0.2}}
\put(4.5,0.3){\circle*{0.1}}
\put(5,0.3){\makebox(0,0)[lc]{\(\ldots\) GI action by \(C_2\)}}
\put(4.5,0){\circle*{0.2}}
\put(5,0){\makebox(0,0)[lc]{\(\ldots\) GI and RI action by \(C_2\)}}
\put(4.5,-0.3){\circle{0.1}}
\put(5,-0.3){\makebox(0,0)[lc]{\(\ldots\) non-\(\sigma\)}}

\put(3.2,0){\makebox(0,0)[lb]{\(\langle 21\rangle\)}}
\put(3.2,-2){\makebox(0,0)[lb]{\(\langle 191\rangle\)}}
\put(3.1,-4){\makebox(0,0)[lb]{\(\langle 1779\rangle\)}}
\multiput(3.2,-6)(0,-2){4}{\makebox(0,0)[lc]{\(1\)}}

\multiput(3,0)(0,-4){4}{\circle{0.2}}
\multiput(3,0)(0,-4){4}{\circle*{0.1}}

\multiput(3,-2)(0,-4){3}{\circle*{0.2}}

\multiput(3,0)(0,-2){6}{\line(0,-1){2}}

\put(3,-12){\vector(0,-1){2}}
\put(3.2,-12.5){\makebox(0,0)[lc]{main line}}

\put(1.2,-2){\makebox(0,0)[lc]{\(\langle 192\rangle\)}}
\put(5.2,-2){\makebox(0,0)[lc]{\(\langle 190\rangle\)}}

\put(1.1,-4){\makebox(0,0)[lc]{\(\langle 1782\rangle\)}}
\put(5.1,-4){\makebox(0,0)[lc]{\(\langle 1780\rangle\)}}
\put(6.1,-4){\makebox(0,0)[lc]{\(\langle 1781\rangle\)}}

\multiput(1.2,-6)(0,-4){2}{\makebox(0,0)[lc]{\(3\)}}
\multiput(1.2,-8)(0,-4){2}{\makebox(0,0)[lc]{\(4\)}}
\multiput(5.2,-6)(0,-2){4}{\makebox(0,0)[lc]{\(2\)}}
\multiput(6.2,-8)(0,-4){2}{\makebox(0,0)[lc]{\(3\)}}

\multiput(1,-2)(0,-4){3}{\circle{0.1}}

\multiput(1,-4)(0,-4){3}{\circle*{0.2}}

\multiput(5,-2)(0,-4){3}{\circle{0.1}}

\multiput(5,-4)(0,-4){3}{\circle{0.2}}
\multiput(5,-4)(0,-4){3}{\circle*{0.1}}
\multiput(6,-4)(0,-4){3}{\circle{0.2}}
\multiput(6,-4)(0,-4){3}{\circle*{0.1}}

\multiput(3,0)(0,-2){6}{\line(-1,-1){2}}
\multiput(3,0)(0,-2){6}{\line(1,-1){2}}
\multiput(3,-2)(0,-4){3}{\line(3,-2){3}}

\multiput(5,-4)(0,-4){2}{\line(1,-1){2}}
\multiput(6,-4)(0,-4){2}{\line(1,-1){2}}

\put(-1,-13){\makebox(0,0)[cc]{pTKT}}
\put(-1,-13.5){\makebox(0,0)[cc]{\(\varkappa\sim\)}}
\put(1,-13){\makebox(0,0)[cc]{\(\mathrm{D}.6\)}}
\put(1,-13.5){\makebox(0,0)[cc]{\((123;1)\)}}
\put(3.5,-13){\makebox(0,0)[cc]{\(\mathrm{e}.14\)}}
\put(3.5,-13.5){\makebox(0,0)[cc]{\((123;0)\)}}
\put(5.5,-13){\makebox(0,0)[cc]{\(\mathrm{E}.12\)}}
\put(5.5,-13.5){\makebox(0,0)[cc]{\((123;4)\)}}

\end{picture}

}

\end{figure}
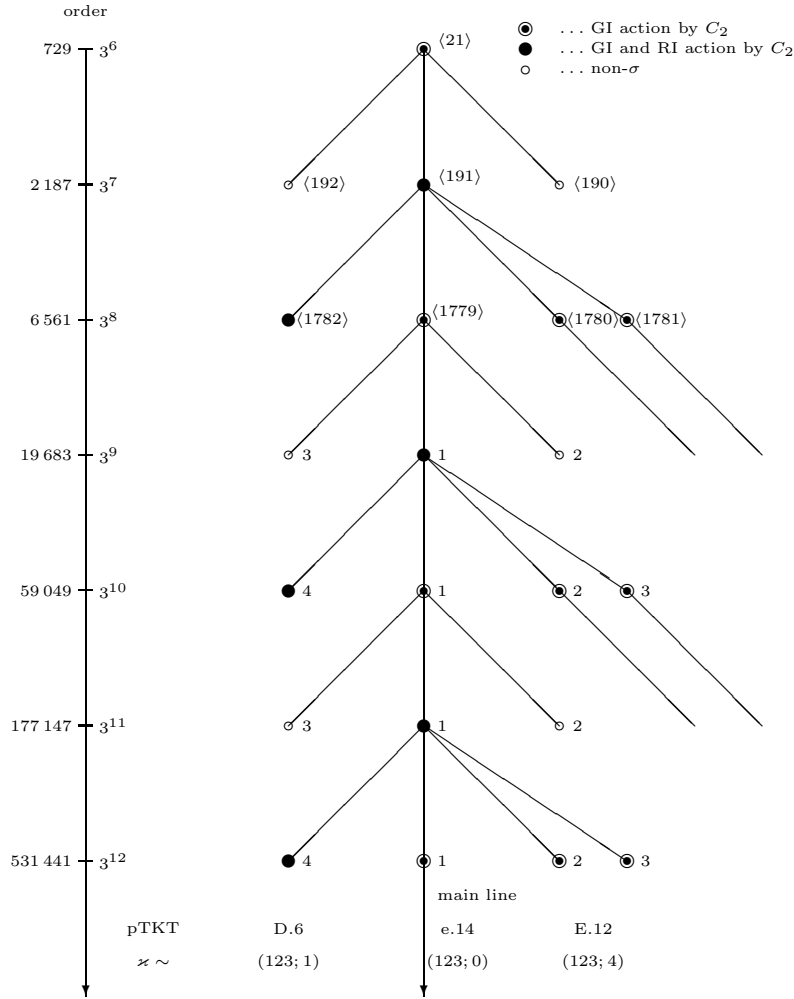


\begin{experiment}
\label{exp:ImagQuadratic}
Among imaginary quadratic fields \(K=\mathbb{Q}(\sqrt{d})\)
with fundamental discriminants \(d<0\)
and \(3\)-class group \(\mathrm{Cl}_3(K)\simeq C_9\times C_3\),
a dominant proportion of
\(\mathbf{46.4}\%\) has a \textbf{metabelian} \(3\)-class field tower with
length \(\ell_3(K)=\mathbf{2}\) and \textbf{Schur} \(\sigma\)-\textbf{group}
\cite{KoVe1975,Ag1998,BBH2017}
\(G=\mathrm{Gal}(\mathrm{F}_3^\infty(K)/K)\simeq\)
either \(\langle 729,14\rangle\) or \(\langle 729,15\rangle\)
having \textit{generator- and relator-inverting} action by \(C_2\) and punctured transfer kernel type \(\mathrm{D}.11\),
\(\varkappa(G)\sim (124;1)\).
Even in the case of equidistribution among the two candidate groups,
the proportion is dominant with \(\mathbf{23.2}\%\) for each.
\end{experiment}


\noindent
We give proofs
for the dominant situation of two two-stage towers (nearly \textit{one half})
and for the significant contribution of eight towers with precisely three stages (more than \(18.8\%\)).

\begin{theorem}
\label{thm:MetabelianTower}
(Two-stage tower.) \\
Let \(K=\mathbb{Q}(\sqrt{d})\) be an imaginary quadratic field
with fundamental discriminant \(d<0\),
\(3\)-class group \(\mathrm{Cl}_3(K)\simeq C_9\times C_3\) and
punctured capitulation type \(\mathrm{D}.11\), \(\varkappa(K)\sim (124;1)\).
\begin{enumerate}
\item
The Galois group \(G_2\) of the second Hilbert \(3\)-class field \(\mathrm{F}_3^2(K)\) of \(K\)
is one of the two unique metabelian Schur \(\sigma\)-groups
\(\mathrm{Gal}(\mathrm{F}_3^2(K)/K)\simeq\langle 729,14\rangle\) or \(\langle 729,15\rangle\) (see Figure
\ref{fig:C9xC3}).
\item
The abelian type invariants of the \(3\)-class groups \(\mathrm{Cl}_3(L_i)\)
of the four unramified cyclic cubic extensions \(L_i\), \(1\le i\le 4\), of \(K\)
are given by \(\alpha(K)\sim (31,31,211;211)\).
\item
The \(3\)-class field tower of \(K\)
stops at the second stage, that is,
\(\mathrm{F}_3^2(K)=\mathrm{F}_3^\infty(K)\) is the maximal unramified pro-\(3\) extension of \(K\).
\end{enumerate}
\end{theorem}

\begin{proof}
The groups \(\langle 729,14\rangle\) or \(\langle 729,15\rangle\)
are unique with Artin pattern \(\alpha(K)=(31,31,211;211)\), \(\varkappa(K)\sim (124;1)\), 
of type \(\mathrm{D}.11\).
\end{proof}


In the next theorem,
a tower with precisely three stages is warranted, independently of the AQI.
For the unambiguous identification of \(G_2\) and \(G_\infty\), however,
a specification of the AQI is mandatory.

\begin{theorem}
\label{thm:ThreeStageTowers}
(Three-stage towers.) \\
Let \(K=\mathbb{Q}(\sqrt{d})\) be an imaginary quadratic field
with fundamental discriminant \(d<0\),
\(3\)-class group \(\mathrm{Cl}_3(K)\simeq C_9\times C_3\) and
punctured capitulation type either
\(\mathrm{D}.10\), \(\varkappa(K)\sim (411;3)\) or
\(\mathrm{C}.4\), \(\varkappa(K)\sim (311;3)\) or
\(\mathrm{D}.5\), \(\varkappa(K)\sim (211;3)\) or
\(\mathrm{D}.6\), \(\varkappa(K)\sim (123;1)\).
\begin{enumerate}
\item
The \(3\)-class field tower of \(K\)
stops at the third stage, that is,
\(\mathrm{F}_3^3(K)=\mathrm{F}_3^\infty(K)\) is the maximal unramified pro-\(3\) extension of \(K\),
for all four assigned pTKTs.
\item
Let the abelian type invariants of the \(3\)-class groups \(\mathrm{Cl}_3(L_i)\)
of the four unramified cyclic cubic extensions \(L_i\), \(1\le i\le 4\), of \(K\)
be denoted by \(\alpha(K)\).
The Galois group \(G_2\) of the second Hilbert \(3\)-class field \(\mathrm{F}_3^2(K)\) of \(K\)
and the Galois group \(G_\infty\) of the Hilbert \(3\)-class field tower \(\mathrm{F}_3^3(K)=\mathrm{F}_3^\infty(K)\)
of \(K\) are given by
\begin{enumerate}
\item
\(G_2\simeq\langle 6561,1683\rangle\) (Figure
\ref{fig:Ord729Id13})
and \(G_\infty\simeq\langle 2187,168\rangle-\#2;2\), \\
if \(\varkappa(K)\sim (311;3)\) and \(\alpha(K)\sim (222,31,31;211)\) (homocyclic type \(\mathrm{C}.4\)),
\item
\(G_2\simeq\langle 6561,1684\rangle\) (Figure
\ref{fig:Ord729Id13})
and \(G_\infty\simeq\langle 2187,168\rangle-\#2;3\), \\
if \(\varkappa(K)\sim (211;3)\) and \(\alpha(K)\sim (222,31,31;211)\) (homocyclic type \(\mathrm{D}.5\)),
\item
\(G_2\simeq\langle 6561,1686\rangle\) (Figure
\ref{fig:Ord729Id13})
and \(G_\infty\simeq\langle 2187,168\rangle-\#2;5\), \\
if \(\varkappa(K)\sim (211;3)\) and \(\alpha(K)\sim (321,31,31;211)\) (heterocyclic type \(\mathrm{D}.5\)),
\item
\(G_2\simeq\langle 6561,1687\rangle\) (Figure
\ref{fig:Ord729Id13})
and \(G_\infty\simeq\langle 2187,168\rangle-\#2;6\), \\
if \(\varkappa(K)\sim (311;3)\) and \(\alpha(K)\sim (321,31,31;211)\) (heterocyclic type \(\mathrm{C}.4\)),
\item
\(G_2\simeq\langle 6561,1689\rangle\) or \(G_2\simeq\langle 6561,1690\rangle\) (Figure
\ref{fig:Ord729Id13}), \\
if \(\varkappa(K)\sim (411;3)\) and \(\alpha(K)\sim (321,31,31;211)\) ( type \(\mathrm{D}.10\)),
\item
\(G_2\simeq\langle 6561,1782\rangle\) (Figure
\ref{fig:Ord729Id21})
or \(G_2\simeq\langle 6561,1744\rangle\), \\
if \(\varkappa(K)\sim (123;1)\) and \(\alpha(K)\sim (31,31,31;321)\) ( type \(\mathrm{D}.6\)).
\end{enumerate}
In each case, \(G_\infty\) is a non-metabelian Schur \(\sigma\)-group with
\(\mathrm{ord}=19\,683\) and \(\mathrm{dl}=3\).
\end{enumerate}
\end{theorem}

\begin{proof}
See Theorem
\ref{thm:Imaginary21}.
\end{proof}


\section{Metabelian Schur \(\sigma\)-groups \(G\) with \(G/G^\prime\simeq (3^e,3)\), \(e\ge 3\)}
\label{s:Metabelian}

\noindent
According to Table
\ref{tbl:ImagQuadratic},
nearly \textit{one half} of the imaginary quadratic fields \(K=\mathbb{Q}(\sqrt{d})\)
with \(3\)-class group \(\mathrm{Cl}_3(K)\simeq (9,3)\)
possesses a \textit{metabelian} \(3\)-class field tower with automorphism group
\(G=\mathrm{Gal}(\mathrm{F}_3^\infty(K)/K)\simeq\mathrm{SmallGroup}(729,i)\), where \(i\in\lbrace 14,15\rbrace\).
More precisely and more generally,
this trend continues for bigger commutator quotients \(G/G^\prime\simeq (3^e,3)\)
and absolute discriminants \(\lvert d\rvert\) below a given upper bound \(B\): \\
\(406\) among \(875\), that is \(46.4\%\), with \(G\simeq\langle 729,14\vert 15\rbrace\), \(G/G^\prime\simeq (9,3)\), when \(B=10^6\), \\
\(433\) among \(930\), that is \(46.56\%\), with \(G\simeq\langle 2187,121\vert 122\rbrace\), \(G/G^\prime\simeq (27,3)\), when \(B=3\cdot 10^6\), \\
\(999\) among \(2174\), that is \(45.95\%\), with \(G\simeq\langle 6561,975\vert 976\rbrace\), \(G/G^\prime\simeq (81,3)\), when \(B=20\cdot 10^6\).

This high proportion of \textit{metabelian} Schur \(\sigma\)-groups
suggested to search for a general theoretical statement concerning metabelian Schur \(\sigma\)-groups \(G\)
(with derived length \(\mathrm{dl}(G)=2\)).
Indeed, we succeeded in finding another \textit{periodicity},
underpinned with two \textit{infinite limit groups} by M. F. Newman,
and justified rigorously by a \textit{parameterized presentation}.

\begin{theorem}
\label{thm:MetabelianSchurSigma}
For all integers \(e\ge 3\),
the \textbf{unique} pair of \textbf{metabelian} Schur \(\sigma\)-groups \(G\)
with commutator quotient \(G/G^\prime\simeq (3^e,3)\hat{=}(e1)\) is given by
the \textbf{periodic sequence}
\begin{equation}
\label{eqn:D11}
\mathrm{SmallGroup}(729,8)(-\#1;1)^{e-3}-\#1;i, \qquad i\in\lbrace 2,3\rbrace.
\end{equation}
\end{theorem}

\begin{corollary}
\label{cor:MetabelianSchurSigma}
The order of these groups is \(\#G=3^{4+e}\),
and their Artin pattern \((\varkappa,\alpha)\) is given by
constant punctured transfer kernel type \(\mathrm{D}.11\), \(\varkappa\sim (124;1)\),
and increasing abelian quotient invariants, \(\alpha\sim ((e+1)1,(e+1)1,e11;e11)\).
\end{corollary}


\begin{theorem}
\label{thm:MetabelianLimit}
For all integers \(e\ge 3\),
the unique pair of metabelian Schur \(\sigma\)-groups \(G\)
with commutator quotient \(G/G^\prime\simeq (3^e,3)\hat{=}(e1)\) is alternatively given by
\(G\simeq (L_{11}/P_e(L_{11}))-\#1;i,\ i\in\lbrace 2,3\rbrace\),
where the \textbf{infinite limit} group \(L_{11}\) is given by the finite presentation
\begin{equation}
\label{eqn:D11Limit}
L_{11}=\langle
a,t,u\mid
\lbrack t,a\rbrack=u,
\lbrack u,a\rbrack=t^3,
t^3=\lbrack u,t\rbrack,
u^3=1
\rangle.
\end{equation}
\end{theorem}


\begin{theorem}
\label{thm:MetabelianQuotient}
For all integers \(e\ge 3\),
the unique pair of metabelian Schur \(\sigma\)-groups \(G\)
with commutator quotient \(G/G^\prime\simeq (3^e,3)\hat{=}(e1)\) is alternatively given by
\(G\simeq L_{16}/\langle a^{\mp 3^e}\cdot t^3\cdot\lbrack t,a,t\rbrack\rangle\)
where the \textbf{infinite limit} group \(L_{16}\) is given by the finite presentation
\begin{equation}
\label{eqn:b16Limit}
L_{16}=\langle
a,t\mid
\lbrack t,a,a\rbrack=t^3,
\lbrack t,a\rbrack^3=1,
\lbrack t,a,t,t\rbrack=1
\rangle.
\end{equation}
Additionally, the \(p\)-parent of the next pair (with \(G/G^\prime\simeq (3^{e+1},3)\)) is given by 
\(L_{16}/\langle t^3\cdot\lbrack t,a,t\rbrack\rangle\).
\end{theorem}


\begin{theorem}
\label{thm:ParametrizedPresentation}
For all integers \(e\ge 3\),
the unique pair of metabelian Schur \(\sigma\)-groups \(G\)
with commutator quotient \(G/G^\prime\simeq (3^e,3)\hat{=}(e1)\) is alternatively given by
the \textbf{parametrized presentation}
\begin{equation}
\label{eqn:ParametrizedPresentation}
G=\langle
x,y,s_2,s_3,t_3,w\mid
s_2=\lbrack y,x\rbrack,
s_3=\lbrack s_2,x\rbrack,
t_3=\lbrack s_2,y\rbrack,
x^{3^e}=w,
y^3=s_3,
R
\rangle,
\end{equation}
where \(R\) denotes either the relation \(t_3=s_3\cdot w\) or \(t_3=s_3\cdot w^2\).
\end{theorem}

\begin{proof}
Theorems
\ref{thm:MetabelianSchurSigma},
\ref{thm:MetabelianLimit} and
\ref{thm:MetabelianQuotient}
are proved by means of Formula
\eqref{eqn:ParametrizedPresentation} in Theorem
\ref{thm:ParametrizedPresentation}.
\end{proof}

\medskip
In particular, for \(e=3\) and \(e=4\), we have the pairs \\
\(\mathrm{SmallGroup}(2187,121)\simeq\langle x,y,s_2,s_3,t_3,w\mid x^{27}=w,y^3=s_3,t_3=s_3\cdot w\rangle\), \\
\(\mathrm{SmallGroup}(2187,122)\simeq\langle x,y,s_2,s_3,t_3,w\mid x^{27}=w,y^3=s_3,t_3=s_3\cdot w^2\rangle\); and \\
\(\mathrm{SmallGroup}(6561,975)\simeq\langle x,y,s_2,s_3,t_3,w\mid x^{81}=w,y^3=s_3,t_3=s_3\cdot w\rangle\), \\
\(\mathrm{SmallGroup}(6561,976)\simeq\langle x,y,s_2,s_3,t_3,w\mid x^{81}=w,y^3=s_3,t_3=s_3\cdot w^2\rangle\).

\medskip
The smallest pair, for \(e=2\), however, is exceptional: \\
\(\mathrm{SmallGroup}(729,14)\simeq\langle x,y,s_2,s_3,t_3\mid x^9=s_3,y^3=t_3\rangle\), \\
\(\mathrm{SmallGroup}(729,15)\simeq\langle x,y,s_2,s_3,t_3\mid x^9=s_3,y^3=t_3^2\rangle\).


\section{Imaginary quadratic fields with bicyclic \(3\)-class group of order \(>27\)}
\label{s:HigherOrder}

\noindent
When we pass from finite \(3\)-groups \(G\)
with simplest non-elementary bicyclic commutator quotient \(G/G^\prime\) of type \((9,3)\)
to situations with order \(>27\),
that is, the cases \((27,3)\), \((81,3)\), \((243,3)\), etc.,
then the construction of the groups by means of the \(p\)-group generation algorithm
\cite{Nm1977,Ob1990}
becomes increasingly difficult,
since the algorithm uses the definition of (parent,descendant)-pairs \((\pi_p(D),D)\) with
\(\pi_p(D)=D/P_k(D)\) in terms of the lower exponent-\(p\) central series \((P_i(D))_{0\le i\le k}\),
whereas the parent \(\pi(D)=D/\gamma_c(D)\) on the coclass tree of \(D\) may have nuclear rank \(n=0\).


In order to emphasize the broad scope of our current investigations,
let us summarize the characteristic properties of all coclass trees
with \textit{metabelian} mainline of type \(\mathrm{d}.10\),
independently of the arbitrary commutator quotient \((9,3)\), \((27,3)\), \((81,3)\), etc.
For fixed commutator quotient this tree is \textit{unique}.
The infinite main line consists of metabelian vertices
with punctured transfer kernel type 
\(\mathrm{d}.10\), \(\varkappa\sim (011;3)\), infinitely capable.
The metabelian vertices of \textit{depth one} with respect to the main line have either \\
type \(\mathrm{B}.2\), \(\varkappa\sim (111;3)\), finitely capable for fixed coclass, infinitely for all step sizes, or \\
type \(\mathrm{D}.5\), \(\varkappa\sim (211;3)\), terminal, or \\
type \(\mathrm{C}.4\), \(\varkappa\sim (311;3)\), terminal, or \\
type \(\mathrm{D}.10\), \(\varkappa\sim (411;3\)), terminal.

Type \(\mathrm{B}.2\) gives rise to brushwood descendants
(if multifurcation is admitted, then
the derived length is probably unbounded),
which prohibits straightforward statements about the tower length
(except that it must be at least equal to three).
However,
we believe that the types \(\mathrm{D}.5\), \(\mathrm{C}.4\) and \(\mathrm{D}.10\)
of imaginary quadratic fields
are always associated with a \textit{tower of precise length three}.

When we assign a fixed punctured transfer target type \(\alpha\)
(i.e. the abelian quotient invariants of maximal subgroups),
then our belief can be proven rigorously for
commutator quotient \((9,3)\)
(Theorem \ref{thm:Imaginary21}, June 2013),
\((27,3)\)
(Theorem \ref{thm:Imaginary31}, 16 July 2021),
\((81,3)\)
(Theorem \ref{thm:Imaginary41}, 20 July 2021),
\((243,3)\)
(Theorem \ref{thm:Imaginary51}, 28 July 2021),
\((729,3)\) and higher
(Theorems \ref{thm:Periodicity} and \ref{thm:Limits}, 30 July 2021).


In the following theorems,
a tower with precisely three stages is warranted, independently of the AQI.
For the unambiguous identification of the metabelianization \(M=G/G^{\prime\prime}\)
and the group \(G\), however,
a specification of the AQI is mandatory.

\begin{theorem}
\label{thm:Imaginary21}
An imaginary quadratic field \(K=\mathbb{Q}(\sqrt{d})\), \(d < 0\),
with \(3\)-class group of type \((9,3)\) and
one of the following six kinds of Artin pattern \(\mathrm{AP}(K)=(\varkappa,\alpha)\), \\
type \(\mathrm{C}.4\), \(\varkappa\sim (311;3)\), \(\alpha\sim (222,31,31;211)\) with homocyclic \(1\)st component or \\
type \(\mathrm{D}.5\), \(\varkappa\sim (211;3)\), \(\alpha\sim (222,31,31;211)\) with homocyclic \(1\)st component or \\
type \(\mathrm{D}.5\), \(\varkappa\sim (211;3)\), \(\alpha\sim (321,31,31;211)\) with heterocyclic \(1\)st component or \\
type \(\mathrm{C}.4\), \(\varkappa\sim (311;3)\), \(\alpha\sim (321,31,31;211)\) with heterocyclic \(1\)st component or \\
type \(\mathrm{D}.10\), \(\varkappa\sim (411;3)\), \(\alpha\sim (321,31,31;211)\) with heterocyclic \(1\)st component or \\
type \(\mathrm{D}.6\), \(\varkappa\sim (123;1)\), \(\alpha\sim (31,31,31;321)\) with heterocyclic \(4\)th component \\
has a \(3\)-class field tower of precise length \(\ell_3(K)=3\).
The tower group \(G=\mathrm{Gal}(\mathrm{F}_3^\infty(K)/K)\) is of order \(3^9\), and
its metabelianization \(M=G/G^{\prime\prime}\) is of order \(3^8\).
\end{theorem}

\begin{proof}
In each case, \(G\) is a non-metabelian Schur \(\sigma\)-group with
\(\mathrm{ord}=19\,683\) and \(\mathrm{dl}=3\).
\begin{enumerate}
\item
\(M\simeq\langle 6561,1683\rangle\) 
and \(G\simeq\langle 2187,168\rangle-\#2;2\), \\
if \(\varkappa(K)\sim (311;3)\) and \(\alpha(K)\sim (222,31,31;211)\) (homocyclic),
\item
\(M\simeq\langle 6561,1684\rangle\) 
and \(G\simeq\langle 2187,168\rangle-\#2;3\), \\
if \(\varkappa(K)\sim (211;3)\) and \(\alpha(K)\sim (222,31,31;211)\) (homocyclic),
\item
\(M\simeq\langle 6561,1686\rangle\) 
and \(G\simeq\langle 2187,168\rangle-\#2;5\), \\
if \(\varkappa(K)\sim (211;3)\) and \(\alpha(K)\sim (321,31,31;211)\) (heterocyclic),
\item
\(M\simeq\langle 6561,1687\rangle\) 
and \(G\simeq\langle 2187,168\rangle-\#2;6\), \\
if \(\varkappa(K)\sim (311;3)\) and \(\alpha(K)\sim (321,31,31;211)\) (heterocyclic),
\item
\(M\simeq\langle 6561,1689\rangle\), \(G\simeq\langle 2187,168\rangle-\#2;8\) or
\(M\simeq\langle 6561,1690\rangle\), \(G\simeq\langle 2187,168\rangle-\#2;9\), \\
if \(\varkappa(K)\sim (411;3)\) and \(\alpha(K)\sim (321,31,31;211)\),
\item
\(M\simeq\langle 6561,1744\rangle\), \(G\simeq\langle 2187,181\rangle-\#2;4\) or
\(M\simeq\langle 6561,1782\rangle\), \(G\simeq\langle 2187,191\rangle-\#2;4\), \\
if \(\varkappa(K)\sim (123;1)\) and \(\alpha(K)\sim (31,31,31;321)\). \qedhere
\end{enumerate}
\end{proof}


\begin{example}
\label{exm:Imaginary21}
Concrete realizations of Artin patterns in Theorem
\ref{thm:Imaginary21}: \\
type \(\mathrm{C}.4\), homocyclic, for \(d = -42\,567\), \\
type \(\mathrm{D}.5\), homocyclic, for \(d = -116\,419\), \\
type \(\mathrm{D}.5\), heterocyclic, for \(d = -17\,723\), \\
type \(\mathrm{C}.4\), heterocyclic, for \(d = -35\,331\), \\
type \(\mathrm{D}.10\) for \(d = -11\,651\), \\
type \(\mathrm{D}.6\) for \(d = -31\,983\).
\end{example}


\begin{theorem}
\label{thm:Imaginary31}
An imaginary quadratic field \(K=\mathbb{Q}(\sqrt{d})\), \(d < 0\),
with \(3\)-class group of type \((27,3)\) and
one of the following four kinds of Artin pattern \(\mathrm{AP}(K)=(\varkappa,\alpha)\), \\
type \(\mathrm{D}.10\), \(\varkappa\sim (411;3)\), \(\alpha\sim (322,41,41;311)\) with homocyclic \(1\)st component or \\
type \(\mathrm{D}.5\), \(\varkappa\sim (211;3)\), \(\alpha\sim (421,41,41;311)\) with heterocyclic \(1\)st component or \\
type \(\mathrm{C}.4\), \(\varkappa\sim (311;3)\), \(\alpha\sim (421,41,41;311)\) with heterocyclic \(1\)st component or \\
type \(\mathrm{D}.6\), \(\varkappa\sim (123;1)\), \(\alpha\sim (41,41,41;322)\) with homocyclic \(4\)th component \\
has a \(3\)-class field tower of precise length \(\ell_3(K)=3\).
The tower group \(G=\mathrm{Gal}(\mathrm{F}_3^\infty(K)/K)\) is of order \(3^{10}\), and
its metabelianization \(M=G/G^{\prime\prime}\) is of order \(3^9\).
\end{theorem}

\begin{proof}
It suffices to give the isomorphism classes of \(G\) and \(M\).
The unique metabelian \(3\)-group with commutator quotient \((27,3)\),
order \(3^8\), and punctured transfer kernel type \(\mathrm{d}.10\), \(\varkappa\sim (011;3)\),
is given by the fork \(B:=\mathrm{SmallGroup}(6561,98)\).
It has nuclear rank \(n(B)=2\) and thus causes a bifurcation
which uniquely determines \(G\) (and thus of course also \(M=G/G^{\prime\prime}\)),
for all given Artin patterns, with two solutions (Schur \(\sigma\)-candidates) each: \\
\(G\simeq B-\#2;2\) resp. \(3\), \(M\simeq B-\#1;4\) resp. \(5\), for \(\varkappa\sim (411;3)\), \(\alpha\sim (322,41,41;311)\), \\
\(G\simeq B-\#2;6\) resp. \(8\), \(M\simeq B-\#1;8\) resp. \(10\), for \(\varkappa\sim (211;3)\), \(\alpha\sim (421,41,41;311)\), \\
\(G\simeq B-\#2;5\) resp. \(9\), \(M\simeq B-\#1;7\) resp. \(11\), for \(\varkappa\sim (311;3)\), \(\alpha\sim (421,41,41;311)\).

There exist two metabelian \(3\)-groups with commutator quotient \((27,3)\),
order \(3^8\), and punctured transfer kernel type \(\mathrm{e}.14\), \(\varkappa\sim (123;0)\),
namely the forks \(B_1:=\mathrm{SmallGroup}(6561,88)\) and \(B_2:=\mathrm{SmallGroup}(6561,91)\).
They have nuclear rank \(n(B_i)=2\).
Their bifurcation determines two Schur \(\sigma\)-candidates for \(G\): \\
\(G\simeq B_i-\#2;4\), and \(M\simeq B_i-\#1;4\) with \(1\le i\le 2\), for \(\varkappa\sim (123;1)\), \(\alpha\sim (41,41,41;322)\).
\end{proof}


\begin{example}
\label{exm:Imaginary31}
Concrete realizations of Artin patterns in Theorem
\ref{thm:Imaginary31}: \\
type \(\mathrm{D}.10\) for \(d = -110\,059\), \\
type \(\mathrm{D}.5\) for \(d = -382\,232\), \\
type \(\mathrm{C}.4\) for \(d = -41\,631\), \\
type \(\mathrm{D}.6\) for \(d = -155\,224\).
\end{example}


\begin{theorem}
\label{thm:Imaginary41}
An imaginary quadratic field \(K=\mathbb{Q}(\sqrt{d})\), \(d < 0\),
with \(3\)-class group of type \((81,3)\) and
one of the following four kinds of Artin pattern \(\mathrm{AP}(K)=(\varkappa,\alpha)\), \\
type \(\mathrm{D}.10\), \(\varkappa\sim (411;3)\), \(\alpha\sim (422,51,51;411)\) with homocyclic \(1\)st component or \\
type \(\mathrm{D}.5\), \(\varkappa\sim (211;3)\), \(\alpha\sim (521,51,51;411)\) with heterocyclic \(1\)st component or \\
type \(\mathrm{C}.4\), \(\varkappa\sim(311;3)\), \(\alpha\sim (521,51,51;411)\) with heterocyclic \(1\)st component or \\
type \(\mathrm{D}.6\), \(\varkappa\sim (123;1)\), \(\alpha\sim (51,51,51;422)\) with homocyclic \(4\)th component \\
has a \(3\)-class field tower of precise length \(\ell_3(K)=3\).
The tower group \(G=\mathrm{Gal}(\mathrm{F}_3^\infty(K)/K)\) is of order \(3^{11}\), and
its metabelianization \(M=G/G^{\prime\prime}\) is of order \(3^{10}\).
\end{theorem}

\begin{proof}
For the leading three Artin patterns,
we start with \(A_1:=\mathrm{SmallGroup}(6561,93)\),
a metabelian \(3\)-group with coclass \(\mathrm{cc}(A_1)=4\) and Artin pattern
\(\varkappa\sim (000;0)\), \(\alpha\sim (421,41,41;311)\).
It has nuclear rank \(n(A_1)=2\).
The metabelianizations \(M\) of the \(3\)-class field tower groups \(G\) are given by \\
\(A_1-\#2;7\) or \(8\) for type \(\mathrm{D}.10\), \(\varkappa\sim (411;3)\), \(\alpha\sim (422,51,51;411)\), \\
\(A_1-\#2;11\) or \(13\) for type \(\mathrm{D}.5\), \(\varkappa\sim (211;3)\), \(\alpha\sim (521,51,51;411)\), \\
\(A_1-\#2;10\) or \(14\) for type \(\mathrm{C}.4\), \(\varkappa\sim(311;3)\), \(\alpha\sim (521,51,51;411)\). \\
They all have nuclear rank \(n(M)=1\)
and a unique terminal descendant \(M-\#1;1\),
which is exactly the Schur \(\sigma\)-group \(G\).

For the trailing Artin pattern,
we start with \(A_2:=\mathrm{SmallGroup}(6561,85)\),
a metabelian \(3\)-group with coclass \(\mathrm{cc}(A_2)=4\) and Artin pattern
\(\varkappa\sim (000;0)\), \(\alpha\sim (41,41,41;322)\).
It has nuclear rank \(n(A_2)=2\).
The metabelianizations \(M\) of the \(3\)-class field tower groups \(G\) are given by \\
\(A_2-\#2;8\) or \(12\) for type \(\mathrm{D}.6\), \(\varkappa\sim (123;1)\), \(\alpha\sim (51,51,51;422)\). \\
They have nuclear rank \(n(M)=1\)
and a unique terminal descendant \(M-\#1;1\),
which is exactly the Schur \(\sigma\)-group \(G\).
\end{proof}


\begin{example}
\label{exm:Imaginary41}
Concrete realizations of Artin patterns in Theorem
\ref{thm:Imaginary41}: \\
type \(\mathrm{D}.10\) for \(d = -469\,283\), \\
type \(\mathrm{D}.5\) for \(d = -584\,411\), \\
type \(\mathrm{C}.4\) for \(d = -617\,363\), \\
type \(\mathrm{D}.6\) for \(d = -548\,939\).
\end{example}


\begin{remark}
\label{rmk:Imaginary41}
The common parent, with respect to the usual lower central,
of the leading \(12\) groups \(G\) and \(M\) is \(B=A_1-\#1;2\) with \(n(B)=1\).
But \(B\) is useless for the construction process,
since it has only \(3\) terminal descendants with Artin pattern of type
\(\mathrm{d}.10\), \(\varkappa\sim (011;3)\), \(\alpha\sim (421,51,51;411)\).
\end{remark}


\begin{theorem}
\label{thm:Imaginary51}
An imaginary quadratic field \(K=\mathbb{Q}(\sqrt{d})\), \(d < 0\),
with \(3\)-class group of type \((243,3)\) and
one of the following four kinds of Artin pattern \(\mathrm{AP}(K)=(\varkappa,\alpha)\), \\
type \(\mathrm{D}.10\), \(\varkappa\sim (411;3)\), \(\alpha\sim (522,61,61;511)\) with homocyclic \(1\)st component or \\
type \(\mathrm{D}.5\), \(\varkappa\sim (211;3)\), \(\alpha\sim (621,61,61;511)\) with heterocyclic \(1\)st component or \\
type \(\mathrm{C}.4\), \(\varkappa\sim(311;3)\), \(\alpha\sim (621,61,61;511)\) with heterocyclic \(1\)st component or \\
type \(\mathrm{D}.6\), \(\varkappa\sim (123;1)\), \(\alpha\sim (61,61,61;522)\) with homocyclic \(4\)th component \\
has a \(3\)-class field tower of precise length \(\ell_3(K)=3\).
The tower group \(G=\mathrm{Gal}(\mathrm{F}_3^\infty(K)/K)\) is of order \(3^{12}\), and
its metabelianization \(M=G/G^{\prime\prime}\) is of order \(3^{11}\).
\end{theorem}

\begin{proof}
First, we start with \(A_2:=\mathrm{SmallGroup}(6561,93)-\#2;2\),
a metabelian \(3\)-group with coclass \(\mathrm{cc}(A_2)=5\) and Artin pattern
\(\varkappa\sim (400;0)\), \(\alpha\sim (522,51,51;411)\).
It has nuclear rank \(n(A_2)=1\).
The metabelianizations \(M\) of the \(3\)-class field tower groups \(G\) are given by \\
\(A_2-\#1;2\) or \(3\) for type \(\mathrm{D}.10\), \(\varkappa\sim (411;3)\), \(\alpha\sim (522,61,61;511)\).

Then, we start with \(A_4:=\mathrm{SmallGroup}(6561,93)-\#2;4\),
a metabelian \(3\)-group with coclass \(\mathrm{cc}(A_4)=5\) and Artin pattern
\(\varkappa\sim (000;0)\), \(\alpha\sim (521,51,51;411)\).
It has nuclear rank \(n(A_4)=1\).
The metabelianizations \(M\) of the \(3\)-class field tower groups \(G\) are given by \\
\(A_4-\#1;2\) or \(3\) for type \(\mathrm{D}.5\), \(\varkappa\sim (211;3)\), \(\alpha\sim (621,61,61;511)\).

Next, we start with \(A_5:=\mathrm{SmallGroup}(6561,93)-\#2;5\),
a metabelian \(3\)-group with coclass \(\mathrm{cc}(A_5)=5\) and Artin pattern
\(\varkappa\sim (000;0)\), \(\alpha\sim (521,51,51;411)\).
It has nuclear rank \(n(A_5)=1\).
The metabelianizations \(M\) of the \(3\)-class field tower groups \(G\) are given by \\
\(A_5-\#1;2\) or \(3\) for type \(\mathrm{C}.4\), \(\varkappa\sim(311;3)\), \(\alpha\sim (621,61,61;511)\).

Finally, we start with \(A_0:=\mathrm{SmallGroup}(6561,85)-\#2;4\),
a metabelian \(3\)-group with coclass \(\mathrm{cc}(A_0)=5\) and Artin pattern
\(\varkappa\sim (000;0)\), \(\alpha\sim (51,51,51;422)\).
It has nuclear rank \(n(A_0)=1\).
The metabelianizations \(M\) of the \(3\)-class field tower groups \(G\) are given by \\
\(A_0-\#1;2\) or \(3\) for type \(\mathrm{D}.6\), \(\varkappa\sim (123;1)\), \(\alpha\sim (61,61,61;522)\).

All these metabelian groups \(M\) have nuclear rank \(n(M)=1\)
and a unique terminal descendant \(M-\#1;1\),
which is exactly the Schur \(\sigma\)-group \(G\).
\end{proof}


\begin{example}
\label{exm:Imaginary51}
Concrete realizations of Artin patterns in Theorem
\ref{thm:Imaginary51}: \\
type \(\mathrm{D}.10\) for \(d = -2\,115\,951\), \\
type \(\mathrm{D}.5\) for \(d = -2\,105\,871\), \\
type \(\mathrm{C}.4\) for \(d = -5\,687\,591\), \\
type \(\mathrm{D}.6\) for \(d = -5\,368\,119\).
\end{example}


\begin{example}
\label{exm:Imaginary61}
Concrete realizations of Artin patterns for \(\mathrm{Cl}_3(\mathbb{Q}(\sqrt{d}))\simeq (729,3)\) in Theorem
\ref{thm:Periodicity}: \\
type \(\mathrm{D}.10\) for \(d = -21\,658\,691\), \\
type \(\mathrm{D}.5\) for \(d = -8\,421\,559\), \\
type \(\mathrm{C}.4\) for \(d = -16\,554\,479\), \\
type \(\mathrm{D}.6\) for \(d = -6\,720\,503\).
\end{example}


\section{Periodicity and limits for \(G/G^\prime\simeq (3^e,3)\), \(e\ge 5\)}
\label{s:PeriodicityAndLimits}

\begin{figure}[ht]
\caption{Schur \(\sigma\)-groups \(G\) with commutator quotient \(G/G^\prime\simeq (3^e,3)\), \(1\le e\le 6\)}
\label{fig:SchurSigma}

{\tiny

\setlength{\unitlength}{0.9cm}
\begin{picture}(14,21.5)(-11,-17.5)

\put(-10,2.5){\makebox(0,0)[cb]{order}}

\put(-10,2){\line(0,-1){15.5}}
\multiput(-10.1,2)(0,-1.5){12}{\line(1,0){0.2}}

\put(-10.2,2){\makebox(0,0)[rc]{\(9\)}}
\put(-9.8,2){\makebox(0,0)[lc]{\(3^2\)}}
\put(-10.2,0.5){\makebox(0,0)[rc]{\(27\)}}
\put(-9.8,0.5){\makebox(0,0)[lc]{\(3^3\)}}
\put(-10.2,-1){\makebox(0,0)[rc]{\(81\)}}
\put(-9.8,-1){\makebox(0,0)[lc]{\(3^4\)}}
\put(-10.2,-2.5){\makebox(0,0)[rc]{\(243\)}}
\put(-9.8,-2.5){\makebox(0,0)[lc]{\(3^5\)}}
\put(-10.2,-4){\makebox(0,0)[rc]{\(729\)}}
\put(-9.8,-4){\makebox(0,0)[lc]{\(3^6\)}}
\put(-10.2,-5.5){\makebox(0,0)[rc]{\(2\,187\)}}
\put(-9.8,-5.5){\makebox(0,0)[lc]{\(3^7\)}}
\put(-10.2,-7){\makebox(0,0)[rc]{\(6\,561\)}}
\put(-9.8,-7){\makebox(0,0)[lc]{\(3^8\)}}
\put(-10.2,-8.5){\makebox(0,0)[rc]{\(19\,683\)}}
\put(-9.8,-8.5){\makebox(0,0)[lc]{\(3^9\)}}
\put(-10.2,-10){\makebox(0,0)[rc]{\(59\,049\)}}
\put(-9.8,-10){\makebox(0,0)[lc]{\(3^{10}\)}}
\put(-10.2,-11.5){\makebox(0,0)[rc]{\(177\,147\)}}
\put(-9.8,-11.5){\makebox(0,0)[lc]{\(3^{11}\)}}
\put(-10.2,-13){\makebox(0,0)[rc]{\(531\,441\)}}
\put(-9.8,-13){\makebox(0,0)[lc]{\(3^{12}\)}}
\put(-10.2,-14.5){\makebox(0,0)[rc]{\(1\,594\,323\)}}
\put(-9.8,-14.5){\makebox(0,0)[lc]{\(3^{13}\)}}

\put(-10,-13.5){\vector(0,-1){2}}

\put(-9.1,1.9){\framebox(0.2,0.2){}}
\put(-9,2){\circle*{0.1}}

\put(-9,0.5){\circle{0.2}}
\put(-9,-2.5){\circle{0.2}}
\put(-9,-4){\circle*{0.2}}
\put(-8.5,-5.5){\circle{0.2}}

\put(-7,-1){\circle{0.2}}
\put(-7,-4){\circle{0.2}}
\put(-7,-5.5){\circle*{0.2}}
\put(-6.5,-7){\circle{0.2}}

\put(-5,-4){\circle{0.2}}
\put(-5,-7){\circle*{0.2}}
\put(-4.5,-8.5){\circle{0.2}}

\put(-3,-7){\circle{0.2}}
\put(-3,-10){\circle{0.2}}

\put(-1,-10){\circle{0.2}}
\put(-1,-11.5){\circle{0.2}}

\put(1,-11.5){\circle{0.2}}
\put(1,-13){\circle{0.2}}

\put(-9.1,-7.1){\framebox(0.2,0.2){}}
\put(-7.1,-8.6){\framebox(0.2,0.2){}}
\put(-5.1,-10.1){\framebox(0.2,0.2){}}
\put(-3.1,-11.6){\framebox(0.2,0.2){}}
\put(-1.1,-13.1){\framebox(0.2,0.2){}}
\put(0.9,-14.6){\framebox(0.2,0.2){}}


\put(-9,2){\line(0,-1){1.5}}
\put(-9,0.5){\line(0,-1){3}}
\put(-9,-2.5){\line(0,-1){1.5}}
\put(-9,-4){\line(0,-1){3}}
\put(-9,-4){\line(1,-3){0.5}}

\put(-9,2){\line(2,-3){2}}
\put(-7,-1){\line(0,-1){3}}
\put(-7,-4){\line(0,-1){1.5}}
\put(-7,-5.5){\line(0,-1){3}}
\put(-7,-5.5){\line(1,-3){0.5}}

\put(-7,-1){\line(2,-3){2}}
\put(-5,-4){\line(0,-1){3}}
\put(-5,-7){\line(0,-1){3}}
\put(-5,-7){\line(1,-3){0.5}}

\put(-5,-4){\line(2,-3){2}}
\put(-3,-7){\line(0,-1){3}}
\put(-3,-10){\line(0,-1){1.5}}

\put(-3,-7){\line(2,-3){2}}
\put(-1,-10){\line(0,-1){1.5}}
\put(-1,-11.5){\line(0,-1){1.5}}

\put(-1,-10){\line(4,-3){2}}
\put(1,-11.5){\line(0,-1){1.5}}
\put(1,-13){\line(0,-1){1.5}}


\put(-8.7,2){\makebox(0,0)[lc]{\(\langle 2\rangle\)}}

\put(-8.7,0.5){\makebox(0,0)[lc]{\(\langle 3\rangle\)}}
\put(-8.7,-2.5){\makebox(0,0)[lc]{\(\langle 8\rangle\)}}
\put(-8.7,-4){\makebox(0,0)[lc]{\(\langle 54\rangle\)}}
\put(-8.2,-5.5){\makebox(0,0)[lc]{\(\langle 304\rangle\)}}

\put(-6.7,-1){\makebox(0,0)[lc]{\(\langle 3\rangle\)}}
\put(-6.7,-4){\makebox(0,0)[lc]{\(\langle 13\rangle\)}}
\put(-6.7,-5.5){\makebox(0,0)[lc]{\(\langle 168\rangle\)}}
\put(-6.2,-7){\makebox(0,0)[lc]{\(\langle 1683\rangle\)}}

\put(-4.7,-4){\makebox(0,0)[lc]{\(\langle 7\rangle\)}}
\put(-4.7,-7){\makebox(0,0)[lc]{\(\langle 98\rangle\)}}
\put(-4.2,-8.5){\makebox(0,0)[lc]{\(1;4\)}}

\put(-2.7,-7){\makebox(0,0)[lc]{\(\langle 93\rangle\)}}
\put(-2.7,-10){\makebox(0,0)[lc]{\(2;7\)}}

\put(-0.7,-10){\makebox(0,0)[lc]{\(2;2\), periodic root}}
\put(-0.7,-11.5){\makebox(0,0)[lc]{\(1;2\)}}

\put(1.3,-11.5){\makebox(0,0)[lc]{\(1;1\)}}
\put(1.3,-13){\makebox(0,0)[lc]{\(1;2\)}}

\put(-8.7,-7){\makebox(0,0)[lc]{\(\langle 622\rangle\)}}
\put(-9,-7.5){\makebox(0,0)[lc]{\(\mathrm{E}.8\)}}
\put(-9,-8){\makebox(0,0)[lc]{\((3,3)\)}}

\put(-6.7,-8.5){\makebox(0,0)[lc]{\(2;2\)}}
\put(-7,-9){\makebox(0,0)[lc]{\(\mathrm{C}.4\)}}
\put(-7,-9.5){\makebox(0,0)[lc]{\((9,3)\)}}

\put(-4.7,-10){\makebox(0,0)[lc]{\(2;2\)}}
\put(-5,-10.5){\makebox(0,0)[lc]{\(\mathrm{D}.10\)}}
\put(-5,-11){\makebox(0,0)[lc]{\((27,3)\)}}

\put(-2.7,-11.5){\makebox(0,0)[lc]{\(1;1\)}}
\put(-3,-12){\makebox(0,0)[lc]{\(\mathrm{D}.10\)}}
\put(-3,-12.5){\makebox(0,0)[lc]{\((81,3)\)}}

\put(-0.7,-13){\makebox(0,0)[lc]{\(1;1\)}}
\put(-1,-13.5){\makebox(0,0)[lc]{\(\mathrm{D}.10\)}}
\put(-1,-14){\makebox(0,0)[lc]{\((243,3)\)}}

\put(1.3,-14.5){\makebox(0,0)[lc]{\(1;1\)}}
\put(1,-15){\makebox(0,0)[lc]{\(\mathrm{D}.10\)}}
\put(1,-15.5){\makebox(0,0)[lc]{\((729,3)\)}}

\put(-3,-17.6){\makebox(0,0)[lc]{Legend:}}

\put(-1.6,-17.7){\framebox(0.2,0.2){}}
\put(-1.5,-17.6){\circle*{0.1}}
\put(-1,-17.6){\makebox(0,0)[lc]{\(\ldots\) abelian}}

\put(-1.5,-18){\circle{0.2}}
\put(-1,-18){\makebox(0,0)[lc]{\(\ldots\) metabelian}}

\put(-1.5,-18.4){\circle*{0.2}}
\put(-1,-18.4){\makebox(0,0)[lc]{\(\ldots\) metabelian with bifurcation}}

\put(-1.6,-18.9){\framebox(0.2,0.2){}}
\put(-1,-18.8){\makebox(0,0)[lc]{\(\ldots\) non-metabelian}}

\end{picture}

}

\end{figure}
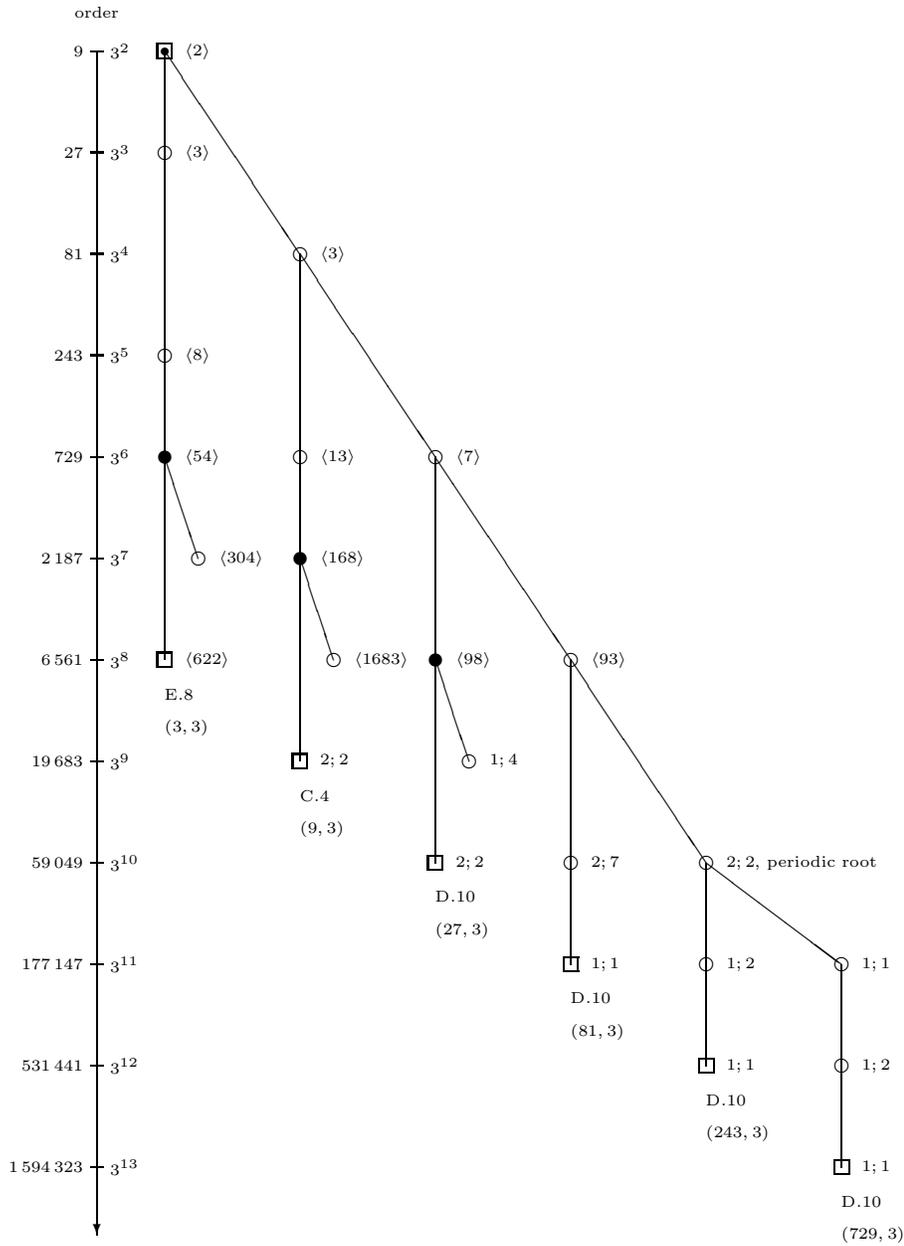

\noindent
In Figure
\ref{fig:SchurSigma},
all directed edges lead from descendants \(D\) to \(p\)-parents \(\pi_p(D)=D/P_{c_p-1}(D)\),
rather than to parents \(\pi(D)=D/\gamma_c(D)\).
The figure admits actual descendant construction.

\noindent
Figure
\ref{fig:SchurSigma}
shows that the construction process for the eight
non-metabelian Schur \(\sigma\)-groups \(G\) with order \(\#G=3^{7+e}\) and
punctured transfer kernel types \(\mathrm{D}.10\), \(\mathrm{C}.4\), \(\mathrm{D}.5\), and \(\mathrm{D}.6\), 
becomes increasingly difficult for the commutator quotients \(G/G^\prime\simeq (27,3)\), \((81,3)\), \((243,3)\).
For the commutator quotient \(G/G^\prime\simeq (729,3)\), however,
an unexpected tranquilization occurs,
and the construction process becomes settled with a simple step size one periodicity.

\newpage

\begin{theorem}
\label{thm:Periodicity}
The four pairs of Schur \(\sigma\)-groups \(G\) with derived length \(\mathrm{dl}(G)=3\),
commutator quotient \(G/G^\prime\simeq (3^e,3)\hat{=}(e1)\), \(e\ge 5\), 
punctured transfer kernel types \(\mathrm{D}.10\), \(\mathrm{D}.5\), \(\mathrm{C}.4\), \(\mathrm{D}.6\),
and order \(\#G=3^{7+e}\)
are given by the following \textbf{bottom up} construction process.
\begin{itemize} 
\item
For type \(\mathrm{D}.10\), \(\varkappa\sim (411;3)\),
\(\alpha\sim (e22,(e+1)1,(e+1)1;e11)\), \\
let \(A_2:=\mathrm{SmallGroup}(6561,93)-\#2;2\), then
\begin{equation}
\label{eqn:D10}
G\simeq A_2(-\#1;1)^{e-5}-\#1;i-\#1;1, \qquad i\in\lbrace 2,3\rbrace.
\end{equation}
\item
For type \(\mathrm{D}.5\), \(\varkappa\sim (211;3)\), 
\(\alpha\sim ((e+1)21,(e+1)1,(e+1)1;e11)\), \\
let \(A_4:=\mathrm{SmallGroup}(6561,93)-\#2;4\), then
\begin{equation}
\label{eqn:D5}
G\simeq A_4(-\#1;1)^{e-5}-\#1;i-\#1;1, \qquad i\in\lbrace 2,3\rbrace.
\end{equation}
\item
For type \(\mathrm{C}.4\), \(\varkappa\sim(311;3)\),
\(\alpha\sim ((e+1)21,(e+1)1,(e+1)1;e11)\), \\
let \(A_5:=\mathrm{SmallGroup}(6561,93)-\#2;5\), then
\begin{equation}
\label{eqn:C4}
G\simeq A_5(-\#1;1)^{e-5}-\#1;i-\#1;1, \qquad i\in\lbrace 2,3\rbrace.
\end{equation}
\item
For type \(\mathrm{D}.6\), \(\varkappa\sim (123;1)\), 
\(\alpha\sim ((e+1)1,(e+1)1,(e+1)1;e22)\), \\
let \(A_0:=\mathrm{SmallGroup}(6561,85)-\#2;4\), then
\begin{equation}
\label{eqn:D6}
G\simeq A_0(-\#1;1)^{e-5}-\#1;i-\#1;1, \qquad i\in\lbrace 2,3\rbrace.
\end{equation}
\end{itemize} 
\end{theorem}


For all types \(\mathrm{D}.10\), \(\mathrm{C}.4\), \(\mathrm{D}.5\), and \(\mathrm{D}.6\),
M. F. Newman has found infinite limit groups
whose quotients give rise to the Schur \(\sigma\)-groups \(G\) of order \(\#G\ge 3^{12}\),
that is, \(e\ge 5\).

\begin{theorem}
\label{thm:Limits}
The four pairs of Schur \(\sigma\)-groups \(G\) with derived length \(\mathrm{dl}(G)=3\),
commutator quotient \(G/G^\prime\simeq (3^e,3)\hat{=}(e1)\), \(e\ge 5\),
punctured transfer kernel types \(\mathrm{D}.10\), \(\mathrm{D}.5\), \(\mathrm{C}.4\), \(\mathrm{D}.6\),
and order \(\#G=3^{7+e}\)
are alternatively given by the following \textbf{top down} construction process.
\begin{itemize} 
\item
For type \(\mathrm{D}.10\), \(\varkappa\sim (411;3)\),
\(\alpha\sim (e22,(e+1)1,(e+1)1;e11)\), \\
let the infinite limit group be given by the finite presentation
\begin{equation}
\label{eqn:D10Limit}
L_{10}=\langle
a,t,u\mid
\lbrack t,a\rbrack=u,
\lbrack u,a\rbrack=\lbrack u,t\rbrack,
\lbrack u,t,u\rbrack=1,
t^3=\lbrack u,t,t,t\rbrack,
u^3=\lbrack u,t,t\rbrack^2\cdot\lbrack u,t,t,t\rbrack
\rangle,
\end{equation}
then
\(G\simeq (L_{10}/P_e(L_{10}))-\#1;i-\#1;1, \qquad i\in\lbrace 2,3\rbrace\).
\item
For type \(\mathrm{D}.5\), \(\varkappa\sim(211;3)\),
\(\alpha\sim ((e+1)21,(e+1)1,(e+1)1;e11)\), \\
let the infinite limit group be given by the finite presentation
\begin{equation}
\label{eqn:D5Limit}
L_5=\langle
a,t,u\mid
\lbrack t,a\rbrack=u,
\lbrack u,a\rbrack=t^3\cdot\lbrack u,t\rbrack,
t^3=\lbrack u,t,t,t\rbrack^{-1}
\rangle,
\end{equation}
then
\(G\simeq (L_5/P_e(L_5))-\#1;i-\#1;1, \qquad i\in\lbrace 2,3\rbrace\).
\item
For type \(\mathrm{C}.4\), \(\varkappa\sim(311;3)\),
\(\alpha\sim ((e+1)21,(e+1)1,(e+1)1;e11)\), \\
let the infinite limit group be given by the finite presentation
\begin{equation}
\label{eqn:C4Limit}
\begin{aligned}
L_4=\langle
a,t,u\mid &
\lbrack t,a\rbrack=u,
\lbrack u,a\rbrack=\lbrack u,t\rbrack\cdot\lbrack u,t,t,t\rbrack^{-1},
\lbrack u,t,u\rbrack=1, \\
& t^3=\lbrack u,t,t,t\rbrack,
u^3=\lbrack u,t,t\rbrack^2\cdot\lbrack u,t,t,t\rbrack,
\lbrack u,t\rbrack^3=\lbrack u,t,t,t\rbrack^2
\rangle,
\end{aligned}
\end{equation}
then
\(G\simeq (L_4/P_e(L_4))-\#1;i-\#1;1, \qquad i\in\lbrace 2,3\rbrace\).
\item
For type \(\mathrm{D}.6\), \(\varkappa\sim(123;1)\),
\(\alpha\sim ((e+1)1,(e+1)1,(e+1)1;e22)\), \\
let the infinite limit group be given by the finite presentation
\begin{equation}
\label{eqn:D6Limit}
L_6=\langle
a,t,u\mid
\lbrack t,a\rbrack=u,
\lbrack u,a\rbrack=t^6\cdot u^6,
\lbrack u,t\rbrack=t^9,
u^9=1,
\lbrack u,t\rbrack^3=1
\rangle,
\end{equation}
then
\(G\simeq (L_6/P_e(L_6))-\#1;i-\#1;1, \qquad i\in\lbrace 2,3\rbrace\).
\end{itemize} 
\end{theorem}

\begin{proof}
The proof consists of the construction of successive descendants of
\(A_2\), \(A_4\), \(A_5\), \(A_0\)
in the way indicated in Theorem
\ref{thm:Periodicity}
by means of the \(p\)-group generation algorithm
\cite{HEO2005}
by Newman
\cite{Nm1977}
and O'Brien
\cite{Ob1990},
which is implemented in the computational algebra system Magma
\cite{BCP1997,BCFS2021,MAGMA2021},
and verifying isomorphism to the descendants of quotients of limit groups
as claimed in Theorem
\ref{thm:Limits}.
\end{proof}


\section{Deterministic laws for CF- and BCF-groups}
\label{s:Laws}

\noindent
The directed edges of the graphs in Figure
\ref{fig:C27xC3}
and
\ref{fig:C81xC3}
are not suitable for the actual construction of the vertices,
since they lead from descendants \(D\) to parents \(P=D/\gamma_c(D)\).
Examplarily, we illustrate some \textit{hidden directed edges}
from descendants \(D\) to \(p\)-parents (ancestors) \(A=D/P_{c_p-1}(D)\).

\begin{theorem}
\label{thm:LowerPCentral}
For each exponent \(e\ge 3\),
we assume that \(a+b=e+1\) with integers \(0\le a-b\le 1\).
Then
\begin{itemize}
\item
there are exactly \textbf{two} \(\mathrm{BCF}\)-groups \(D\)
with punctured transfer kernel type \(\mathrm{e}.14\) and Artin pattern
\(\varkappa\sim (123;0)\), \(\alpha\sim ((e+1)1,(e+1)1,(e+1)1;ab1)\);
their common \(p\)-parent is the \(\mathrm{CF}\)-group
\(A=D/P_e(D)\) with type \(\mathrm{a}.1\) and Artin pattern
\(\varkappa\sim (000;0)\), \(\alpha\sim (e1,e1,e1;ab1)\);
\item
there is a \textbf{unique} \(\mathrm{BCF}\)-group \(D\)
with punctured transfer kernel type \(\mathrm{d}.10\) and Artin pattern
\(\varkappa\sim (110;2)\), \(\alpha\sim ((e+1)1,(e+1)1,e11;e11)\);
its \(p\)-parent is the \(\mathrm{CF}\)-group
\(A=D/P_e(D)\) with type \(\mathrm{a}.1\) and Artin pattern
\(\varkappa\sim (000;0)\), \(\alpha\sim (e1,e1,e11;(e-1)11)\); 
\item
there are exactly \textbf{two} \(\mathrm{BCF}\)-groups \(D\)
with punctured transfer kernel type \(\mathrm{D}.11\) and Artin pattern
\(\varkappa\sim (124;2)\), \(\alpha\sim ((e+1)1,(e+1)1,e11;e11)\);
their common \(p\)-parent is the \(\mathrm{CF}\)-group
\(A=D/P_e(D)\) with type \(\mathrm{b}.16\), Artin pattern
\(\varkappa\sim (004;0)\), \(\alpha\sim (e1,e1,e11;(e-1)11)\).
\end{itemize}
The \(p\)-class of the \(\mathrm{BCF}\)-groups \(D\) is always
\(c_p(D)=e+1\) and their order is \(\#D=3^{4+e}\).
The commutator quotient of all groups \(D\) and \(A\) is
\((3^e,3)\).
\end{theorem}

The groups of the last item are the metabelian Schur \(\sigma\)-groups
in \S\
\ref{s:Metabelian}.

\begin{proof}
We generally denote some crucial commutators by
\(s_2=\lbrack y,x\rbrack\), \(s_3=\lbrack s_2,x\rbrack\), \(t_3=\lbrack s_2,y\rbrack\).
\begin{itemize}
\item
\(D=\langle x,y\mid x^{3^e}=w,y^3=s_3^2,t_3=w \text{ resp. } w^2\rangle\),
\(P_e(D)=\langle w\rangle\simeq C_3\), and thus \\
\(A=D/P_e(D)=\langle x,y\mid x^{3^e}=1,y^3=s_3^2,t_3=1\rangle\).
\item
\(D=\langle x,y\mid x^{3^e}=w,y^3=1,t_3=s_3w\rangle\),
\(P_e(D)=\langle w\rangle\simeq C_3\), and thus \\
\(A=D/P_e(D)=\langle x,y\mid x^{3^e}=1,y^3=1,s_3=t_3\rangle\).
\item
\(D=\langle x,y\mid x^{3^e}=w,y^3=s_3,t_3=s_3w \text{ resp. } s_3w^2\rangle\),
\(P_e(D)=\langle w\rangle\simeq C_3\), and thus \\
\(A=D/P_e(D)=\langle x,y\mid x^{3^e}=1,y^3=s_3=t_3\rangle\). \qedhere
\end{itemize} 
\end{proof}

In particular, we have some relations between groups with SmallGroup identifiers for \(e\in\lbrace 3,4\rbrace\).

\begin{corollary}
\label{cor:LowerPCentralLO7}
The \(p\)-class of the \(\mathrm{BCF}\)-groups
\(D\simeq\langle 2187,i\rangle\) with \(i\in\lbrace 103,104,112,121,122\rbrace\)
and punctured transfer kernel types \(\mathrm{e}.14\), \(\mathrm{d}.10\), \(\mathrm{D}.11\)
is \(c_p=4\), and their \(p\)-parents are the \(\mathrm{CF}\)-groups
\(A=D/P_3(D)\simeq\langle 729,j\rangle\) with \(j\in\lbrace 6,7,8\rbrace\)
and types \(\mathrm{a}.1\), \(\mathrm{b}.16\).
\end{corollary}

\begin{proof}
By specialization of Theorem
\ref{thm:LowerPCentral}
to \(e=3\), we obtain:
\begin{itemize}
\item
For \(D=\langle 2187,103\vert 104\rangle\):
\(P_3(D)\simeq C_3\), and thus
\(A=D/P_3(D)=\langle 729,6\rangle\).
\item
For \(D=\langle 2187,112\rangle\):
\(P_3(D)\simeq C_3\), and thus
\(A=D/P_3(D)=\langle 729,7\rangle\).
\item
For \(D=\langle 2187,121\vert 122\rangle\):
\(P_3(D)\simeq C_3\), and thus
\(A=D/P_3(D)=\langle 729,8\rangle\). \qedhere
\end{itemize} 
\end{proof}

\begin{corollary}
\label{cor:LowerPCentralLO8}
The \(p\)-class of the \(\mathrm{BCF}\)-groups
\(D\simeq\langle 6561,i\rangle\) with \(i\in\lbrace 933,934,953,975,076\rbrace\)
and punctured transfer kernel types \(\mathrm{e}.14\), \(\mathrm{d}.10\), \(\mathrm{D}.11\)
is \(c_p=5\), and their \(p\)-parents are the \(\mathrm{CF}\)-groups
\(A=D/P_4(D)\simeq\langle 2187,j\rangle\) with \(j\in\lbrace 102,111,120\rbrace\)
and types \(\mathrm{a}.1\), \(\mathrm{b}.16\).
\end{corollary}

\begin{proof}
By specialization of Theorem
\ref{thm:LowerPCentral}
to \(e=4\), we get:
\begin{itemize}
\item
For \(D=\langle 6561,933\vert 934\rangle\):
\(P_4(D)\simeq C_3\), and thus
\(A=D/P_4(D)=\langle 2187,102\rangle\).
\item
For \(D=\langle 6561,953\rangle\):
\(P_4(D)\simeq C_3\), and thus
\(A=D/P_4(D)=\langle 2187,111\rangle\).
\item
For \(D=\langle 6561,975\vert 976\rangle\):
\(P_4(D)\simeq C_3\), and thus
\(A=D/P_4(D)=\langle 2187,120\rangle\). \qedhere
\end{itemize} 
\end{proof}

There exist many other similar \textit{deterministic laws} for groups
which were not in the focus of the present paper,
in contrast to those of Theorem
\ref{thm:LowerPCentral}. See the following appendix.


\section{Group theoretic appendix}
\label{s:Appendix}

\noindent
Arithmetical evaluation of the groups investigated in this appendix is difficult.
The exposition is purely group theoretical.
It supplements further interesting \textit{deterministic laws}
for CF-groups and BCF-groups (with moderate, resp. elevated, rank distribution).

\begin{definition}
\label{dfn:Ranks}
By the \textit{rank distribution}
of a pro-\(3\) group \(G\) with bicyclic commutator quotient \(G/G^\prime\)
we understand the punctured quartet \(\varrho(G)\sim (\mathrm{rank}_3(H_i/H_i^\prime))_{(G:H_i)=3}\).
\end{definition}

\begin{theorem}
\label{thm:ModerateBCF}
For each exponent \(e\ge 3\),
we assume that \(a+b=e+1\) with integers \(0\le a-b\le 1\).
Then
\begin{itemize}
\item
there are exactly \textbf{two} \(\mathrm{BCF}\)-groups \(D\)
with punctured transfer kernel type \(\mathrm{B}.7\) and Artin pattern
\(\varkappa\sim (111;4)\), \(\alpha\sim ((e+1)1,(e+1)1,(e+1)1;(e-1)111)\);
their common \(p\)-parent is the \(\mathrm{CF}\)-group
\(A=D/P_e(D)\) with type \(\mathrm{b}.15\) and Artin pattern
\(\varkappa\sim (000;4)\), \(\alpha\sim (e1,e1,e1;(e-1)111)\); 
\item
there are exactly \textbf{two} \(\mathrm{BCF}\)-groups \(D\)
with punctured transfer kernel type \(\mathrm{E}.12\) and Artin pattern
\(\varkappa\sim (123;4)\), \(\alpha\sim ((e+1)1,(e+1)1,(e+1)1;ab1)\);
their common \(p\)-parent is the \(\mathrm{CF}\)-group
\(A=D/P_e(D)\) with type \(\mathrm{b}.15\), Artin pattern
\(\varkappa\sim (000;4)\), \(\alpha\sim (e1,e1,e1;ab1)\).
\end{itemize}
These two cases together with the cases in Theorem
\ref{thm:LowerPCentral}
are summarized in Table
\ref{tbl:ModerateBCF}.
The \(p\)-class of the \(\mathrm{BCF}\)-groups \(D\) is always
\(c_p(D)=e+1\) and their order is \(\#D=3^{4+e}\).
The commutator quotient of all groups \(D\) and \(A\) is
\((3^e,3)\),
and they are tied together by the rank distribution \(\varrho\).
\end{theorem}


\renewcommand{\arraystretch}{1.0}
\begin{table}[ht]
\caption{BCF-groups with moderate rank distribution \(\varrho\)}
\label{tbl:ModerateBCF}
\begin{center}
\begin{tabular}{|c|l|c|c||c||l|c|c|}
\hline
 \multicolumn{4}{|c||}{ BCF-group \(D\) }                                       &             & \multicolumn{3}{|c|}{ CF-group \(A=D/P_e(D)\) }           \\
 \(\#\) & pTKT              & \(\varkappa\) & \(\alpha\)                        & \(\varrho\) & pTKT              & \(\varkappa\) & \(\alpha\)            \\
\hline
 \(2\)  & \(\mathrm{B}.7\)  & \(111;4\)     & \((e+1)1,(e+1)1,(e+1)1;(e-1)111\) & \(2,2,2;4\) & \(\mathrm{b}.15\) & \(000;4\)     & \(e1,e1,e1;(e-1)111\) \\
 \(2\)  & \(\mathrm{E}.12\) & \(123;4\)     & \((e+1)1,(e+1)1,(e+1)1;ab1\)      & \(2,2,2;3\) & \(\mathrm{b}.15\) & \(000;4\)     & \(e1,e1,e1;ab1\)      \\
 \(2\)  & \(\mathrm{e}.14\) & \(123;0\)     & \((e+1)1,(e+1)1,(e+1)1;ab1\)      & \(2,2,2;3\) & \(\mathrm{a}.1\)  & \(000;0\)     & \(e1,e1,e1;ab1\)      \\
 \(1\)  & \(\mathrm{d}.10\) & \(110;2\)     & \((e+1)1,(e+1)1,e11;e11\)         & \(2,2,3;3\) & \(\mathrm{a}.1\)  & \(000;0\)     & \(e1,e1,e11;(e-1)11\) \\
 \(2\)  & \(\mathrm{D}.11\) & \(124;2\)     & \((e+1)1,(e+1)1,e11;e11\)         & \(2,2,3;3\) & \(\mathrm{b}.16\) & \(004;0\)     & \(e1,e1,e11;(e-1)11\) \\
\hline
\end{tabular}
\end{center}
\end{table}


\begin{proof}
We only have to justify the statements in the first two rows.
Everything else has been proved in Theorem
\ref{thm:LowerPCentral}.
\begin{itemize}
\item
\(D=\langle x,y\mid x^{3^e}=w,y^3=1,t_3=w \text{ resp. } w^2\rangle\),
\(P_e(D)=\langle w\rangle\simeq C_3\), and thus \\
\(A=D/P_e(D)=\langle x,y\mid x^{3^e}=1,y^3=1,t_3=1\rangle\).
\item
\(D=\langle x,y\mid x^{3^e}=w,y^3=s_3,t_3=w \text{ resp. } w^2\rangle\),
\(P_e(D)=\langle w\rangle\simeq C_3\), and thus \\
\(A=D/P_e(D)=\langle x,y\mid x^{3^e}=1,y^3=s_3,t_3=1\rangle\). \qedhere
\end{itemize} 
\end{proof}

\begin{example}
By specialization of Theorem
\ref{thm:ModerateBCF}
to \(e=3\), we obtain:
\begin{itemize}
\item
For \(D=\langle 2187,84\vert 85\rangle\):
\(P_3(D)\simeq C_3\), and thus
\(A=D/P_3(D)=\langle 729,4\rangle\).
\item
For \(D=\langle 2187,94\vert 95\rangle\):
\(P_3(D)\simeq C_3\), and thus
\(A=D/P_3(D)=\langle 729,5\rangle\). \qedhere
\end{itemize} 

\noindent
By specialization of Theorem
\ref{thm:ModerateBCF}
to \(e=4\), we get:
\begin{itemize}
\item
For \(D=\langle 6561,876\vert 877\rangle\):
\(P_4(D)\simeq C_3\), and thus
\(A=D/P_4(D)=\langle 2187,83\rangle\).
\item
For \(D=\langle 6561,917\vert 918\rangle\):
\(P_4(D)\simeq C_3\), and thus
\(A=D/P_4(D)=\langle 2187,93\rangle\). \qedhere
\end{itemize} 
\end{example}

We call the construction method of BCF-groups with moderate rank distribution
an \textit{intra-genetic propagation}
since the descendant \(D\) and the \(p\)-parent \(A\) share a common commutator quotient.


Now we come to the \textit{extra-genetic propagation}
of CF-groups and BCF-groups with elevated rank distribution.

\begin{theorem}
\label{thm:CF}
For each exponent \(e\ge 3\),
we assume that \(a+b=e+1\) with integers \(0\le a-b\le 1\)
and \(c+d=e\) with integers \(0\le c-d\le 1\).
Then Table
\ref{tbl:CF}
shows five \(\mathrm{CF}\)-groups \(D\) with commutator quotient \((3^{e+1},3)\)
and their \(p\)-parents \(A\) with commutator quotient \((3^e,3)\).
The \(p\)-class of the \(\mathrm{CF}\)-groups \(D\) is always
\(c_p(D)=e+1\) and their order is \(\#D=3^{4+e}\).
\(D\) and \(A\) are tied together by
the punctured transfer kernel type \(\varkappa\)
and the rank distribution \(\varrho\),
but the second component \(\alpha\) of the Artin pattern \((\varkappa,\alpha)\) is distinct.
\end{theorem}


\renewcommand{\arraystretch}{1.0}
\begin{table}[ht]
\caption{CF-groups and their rank distribution \(\varrho\)}
\label{tbl:CF}
\begin{center}
\begin{tabular}{|c||c|l|c||c|}
\hline
 CF-group \(D\)                &             &                   &               & CF-group \(A=D/P_e(D)\) \\
 \(\alpha\)                    & \(\varrho\) & pTKT              & \(\varkappa\) & \(\alpha\)            \\
\hline
 \((e+1)1,(e+1)1,(e+1)1;e111\) & \(2,2,2;4\) & \(\mathrm{b}.15\) & \(000;4\)     & \(e1,e1,e1;(e-1)111\) \\
 \((e+1)1,(e+1)1,(e+1)1;ab1\)  & \(2,2,2;3\) & \(\mathrm{b}.15\) & \(000;4\)     & \(e1,e1,e1;cd1\)      \\
 \((e+1)1,(e+1)1,(e+1)1;ab1\)  & \(2,2,2;3\) & \(\mathrm{a}.1\)  & \(000;0\)     & \(e1,e1,e1;cd1\)      \\
 \((e+1)1,(e+1)1,(e+1)11;e11\) & \(2,2,3;3\) & \(\mathrm{a}.1\)  & \(000;0\)     & \(e1,e1,e11;(e-1)11\) \\
 \((e+1)1,(e+1)1,(e+1)11;e11\) & \(2,2,3;3\) & \(\mathrm{b}.16\) & \(004;0\)     & \(e1,e1,e11;(e-1)11\) \\
\hline
\end{tabular}
\end{center}
\end{table}


\begin{proof}
As before, we denote the main commutators by
\(s_2=\lbrack y,x\rbrack\), \(s_3=\lbrack s_2,x\rbrack\), \(t_3=\lbrack s_2,y\rbrack\).
\begin{itemize}
\item
\(D=\langle x,y\mid x^{3^{e+1}}=1,x^{3^e}=w,y^3=1,t_3=1\rangle\),
\(P_e(D)=\langle w\rangle\simeq C_3\), and thus \\
\(A=D/P_e(D)=\langle x,y\mid x^{3^e}=1,y^3=1,t_3=1\rangle\).
\item
\(D=\langle x,y\mid x^{3^{e+1}}=1,x^{3^e}=w,y^3=s_3,t_3=1\rangle\),
\(P_e(D)=\langle w\rangle\simeq C_3\), and thus \\
\(A=D/P_e(D)=\langle x,y\mid x^{3^e}=1,y^3=s_3,t_3=1\rangle\).
\item
\(D=\langle x,y\mid x^{3^{e+1}}=1,x^{3^e}=w,y^3=s_3^2,t_3=1\rangle\),
\(P_e(D)=\langle w\rangle\simeq C_3\), and thus \\
\(A=D/P_e(D)=\langle x,y\mid x^{3^e}=1,y^3=s_3^2,t_3=1\rangle\).
\item
\(D=\langle x,y\mid x^{3^{e+1}}=1,x^{3^e}=w,y^3=1,t_3=s_3\rangle\),
\(P_e(D)=\langle w\rangle\simeq C_3\), and thus \\
\(A=D/P_e(D)=\langle x,y\mid x^{3^e}=1,y^3=1,t_3=s_3\rangle\).
\item
\(D=\langle x,y\mid x^{3^{e+1}}=1,x^{3^e}=w,y^3=s_3,t_3=s_3\rangle\),
\(P_e(D)=\langle w\rangle\simeq C_3\), and thus \\
\(A=D/P_e(D)=\langle x,y\mid x^{3^e}=1,y^3=s_3,t_3=s_3\rangle\). \qedhere
\end{itemize} 
\end{proof}

\begin{example}
By specialization of Theorem
\ref{thm:CF}
to \(e=3\), we obtain:
\begin{itemize}
\item
For \(D=\langle 2187,83\rangle\):
\(P_3(D)\simeq C_3\), and thus
\(A=D/P_3(D)=\langle 729,4\rangle\).
\item
For \(D=\langle 2187,93\rangle\):
\(P_3(D)\simeq C_3\), and thus
\(A=D/P_3(D)=\langle 729,5\rangle\). 
\item
For \(D=\langle 2187,102\rangle\):
\(P_3(D)\simeq C_3\), and thus
\(A=D/P_3(D)=\langle 729,6\rangle\). 
\item
For \(D=\langle 2187,111\rangle\):
\(P_3(D)\simeq C_3\), and thus
\(A=D/P_3(D)=\langle 729,7\rangle\). 
\item
For \(D=\langle 2187,120\rangle\):
\(P_3(D)\simeq C_3\), and thus
\(A=D/P_3(D)=\langle 729,8\rangle\). \qedhere
\end{itemize} 
\end{example}


\begin{theorem}
\label{thm:BCF}
For each exponent \(e\ge 3\),
we assume that \(a+b=e+2\) with integers \(0\le a-b\le 1\)
and \(c+d=e+1\) with integers \(0\le c-d\le 1\).
Then Table
\ref{tbl:BCF}
shows four \(\mathrm{BCF}\)-groups \(D\) with commutator quotient \((3^{e+1},3)\)
and their \(p\)-parents \(A\) with commutator quotient \((3^e,3)\).
The \(p\)-class of the \(\mathrm{BCF}\)-groups \(D\) is always
\(c_p(D)=e+1\) and their order is \(\#D=3^{5+e}\).
\(D\) and \(A\) are tied together by
the punctured transfer kernel type \(\varkappa\)
and the rank distribution \(\varrho\),
but the second component \(\alpha\) of the Artin pattern \((\varkappa,\alpha)\) is distinct.
\end{theorem}


\renewcommand{\arraystretch}{1.0}
\begin{table}[ht]
\caption{BCF-groups with elevated rank distribution \(\varrho\)}
\label{tbl:BCF}
\begin{center}
\begin{tabular}{|c||c|l|c||c|}
\hline
 BCF-group \(D\)                  &             &                   &               & BCF-group \(A=D/P_e(D)\) \\
 \(\alpha\)                       & \(\varrho\) & pTKT              & \(\varkappa\) & \(\alpha\)               \\
\hline
 \((e+1)11,(e+1)11,(e+1)11;e111\) & \(3,3,3;4\) & \(\mathrm{b}.15\) & \(000;4\)     & \(e11,e11,e11;(e-1)111\) \\
 \((e+1)11,(e+1)11,(e+1)11;ab1\)  & \(3,3,3;3\) & \(\mathrm{b}.31\) & \(044;4\)     & \(e11,e11,e11;cd1\)      \\
 \((e+1)11,(e+1)11,(e+1)11;ab1\)  & \(3,3,3;3\) & \(\mathrm{c}.27\) & \(044;0\)     & \(e11,e11,e11;cd1\)      \\
 \((e+1)11,(e+1)11,(e+1)11;e111\) & \(3,3,3;4\) & \(\mathrm{A}.20\) & \(444;4\)     & \(e11,e11,e11;(e-1)111\) \\
\hline
\end{tabular}
\end{center}
\end{table}


\begin{proof}
Again, we put
\(s_2=\lbrack y,x\rbrack\), \(s_3=\lbrack s_2,x\rbrack\), \(t_3=\lbrack s_2,y\rbrack\).
\begin{itemize}
\item
\(D=\langle x,y\mid x^{3^{e+1}}=1,x^{3^e}=w,y^3=1\rangle\),
\(P_e(D)=\langle w\rangle\simeq C_3\), and thus \\
\(A=D/P_e(D)=\langle x,y\mid x^{3^e}=1,y^3=1\rangle\).
\item
\(D=\langle x,y\mid x^{3^{e+1}}=1,x^{3^e}=w,y^3=s_3\rangle\),
\(P_e(D)=\langle w\rangle\simeq C_3\), and thus \\
\(A=D/P_e(D)=\langle x,y\mid x^{3^e}=1,y^3=s_3\rangle\).
\item
\(D=\langle x,y\mid x^{3^{e+1}}=1,x^{3^e}=w,y^3=s_3^2\rangle\),
\(P_e(D)=\langle w\rangle\simeq C_3\), and thus \\
\(A=D/P_e(D)=\langle x,y\mid x^{3^e}=1,y^3=s_3^2\rangle\).
\item
\(D=\langle x,y\mid x^{3^{e+1}}=1,x^{3^e}=w,y^3=t_3\rangle\),
\(P_e(D)=\langle w\rangle\simeq C_3\), and thus \\
\(A=D/P_e(D)=\langle x,y\mid x^{3^e}=1,y^3=t_3\rangle\). \qedhere
\end{itemize} 
\end{proof}

\begin{example}
By specialization of Theorem
\ref{thm:BCF}
to \(e=3\), we obtain:
\begin{itemize}
\item
For \(D=\langle 6561,200\rangle\):
\(P_3(D)\simeq C_3\), and thus
\(A=D/P_3(D)=\langle 2187,2\rangle\).
\item
For \(D=\langle 6561,216\rangle\):
\(P_3(D)\simeq C_3\), and thus
\(A=D/P_3(D)=\langle 2187,3\rangle\). 
\item
For \(D=\langle 6561,229\rangle\):
\(P_3(D)\simeq C_3\), and thus
\(A=D/P_3(D)=\langle 2187,4\rangle\). 
\item
For \(D=\langle 6561,242\rangle\):
\(P_3(D)\simeq C_3\), and thus
\(A=D/P_3(D)=\langle 2187,5\rangle\). \qedhere
\end{itemize} 
\end{example}


\section{Conclusion}
\label{s:Conclusion}

\noindent
Since the relevant \(3\)-groups \(G\) in this paper
are realized as Galois groups \(G\simeq\mathrm{Gal}(\mathrm{F}_3^k(K)/K)\)
of iterated unramified Hilbert \(3\)-class fields (\(2\le k\le\infty\))
of \textit{imaginary} quadratic fields \(K\),
they must be \(\sigma\)-groups.
Therefore, only \textit{every other} branch of a descendant tree
consists of \textit{admissible} vertices.
We say the groups on the first admissible branch are in the \textit{ground state}
and the groups on higher admissible branches  are in \textit{excited states}.

Using this terminology, we can easily point out the entirely different nature
of the infinite limit groups in the articles
\cite{ELNO2013,Ma2018}
and those in \S\S\
\ref{s:Metabelian}
and
\ref{s:PeriodicityAndLimits}
of the present paper.
\begin{itemize}
\item
The limit groups in
\cite{ELNO2013} and
\cite{Ma2018}
admit the construction of \textit{all excited states} with \textit{fixed} commutator quotient
(\((3,3)\) or \((9,3)\) in
\cite[Thm. 4.1]{ELNO2013},
only \((3,3)\) in
\cite[Thm. 3.5]{Ma2018}).
\item
The limit groups in the present paper
admit the construction of \textit{the ground state} with \textit{varying} commutator quotient
(\((3^e,3)\) with \(e\ge 5\) in Theorem
\ref{thm:Limits}).
\end{itemize}

Our theory in \S\
\ref{s:Metabelian}
is complete, since groups with punctured transfer kernel type \(\mathrm{D}.11\)
can be called \textit{sporadic} (outside of coclass trees)
and exist only in the ground state,
namely as \textit{metabelian} Schur \(\sigma\)-groups \(G\).
Since these groups are of class \(\mathrm{cl}(G)=3\),
the periodicity in Theorems
\ref{thm:MetabelianSchurSigma},
\ref{thm:MetabelianLimit},
\ref{thm:MetabelianQuotient}
and
\ref{thm:ParametrizedPresentation}
sets in with \(e=3\).
With the aid of Theorem
\ref{thm:LowerPCentral},
the metabelian Schur \(\sigma\)-groups
are embedded into more general deterministic laws,
which express a remarkable relationship
between BCF groups (with moderate rank distribution) and CF-groups.

The theory in \S\
\ref{s:PeriodicityAndLimits}
is not complete,
since Theorems
\ref{thm:Periodicity}
and
\ref{thm:Limits}
only deal with the \textit{ground state} of
groups with punctured transfer kernel types
\(\mathrm{D}.10\), \(\mathrm{C}.4\), \(\mathrm{D}.5\), \(\mathrm{D}.6\), 
which are of class \(\mathrm{cl}(G)=5\),
whence periodicity sets in with \(e=5\).
But these groups are \textit{periodic} vertices of coclass trees.

In theorems concerning the \textit{first excited state} of
groups with punctured transfer kernel types
\(\mathrm{D}.10\), \(\mathrm{C}.4\), \(\mathrm{D}.5\), \(\mathrm{D}.6\), 
which are of class \(\mathrm{cl}(G)=7\),
periodicity would set in with \(e=7\), and so on.

Figure
\ref{fig:SchurSigma}
illuminates the broad range of \textit{non-metabelian} Schur \(\sigma\)-groups \(G\)
with growing commutator quotients \(G/G^\prime\),
beginning with \((3,3)\) in
\cite{Ma2018}
and \((9,3)\) in Theorem
\ref{thm:Imaginary21},
where descendants are still constructed with (usual) parents in coclass trees,
over the increasingly irregular cases \((27,3)\) in Theorem
\ref{thm:Imaginary31}
and \((81,3)\) in Theorem
\ref{thm:Imaginary41},
up to the \textit{new periodicity} for \((3^e,3)\) with \(e\ge 5\),
which starts in Theorem
\ref{thm:Imaginary51}
and continues in Theorem
\ref{thm:Periodicity},
where descendants are constructed with \(p\)-parents,
the coclass trees begin to hide,
and the bifurcations degenerate to simple descendant relations.


\section{Acknowledgement}
\label{s:Thanks}

\noindent
The author thanks Professor M. F. Newman
from the Australian National University in Canberra, Australian Capital Territory,
for the infinite limit groups in \S\S\
\ref{s:Metabelian}
and
\ref{s:PeriodicityAndLimits}.



\end{document}